\documentclass[
               DIV=15,
               ]{scrartcl}  

\usepackage{color} 
\usepackage{psfrag} 
\usepackage{subcaption}
\usepackage{hyperref}
\usepackage{amsmath}
\usepackage{amssymb}
\usepackage{amsthm}
\usepackage{graphics}
\usepackage{todonotes}
\usepackage{amsfonts}  
\usepackage{adjustbox}    
\usepackage{mathtools}
\usepackage{xcolor}        
\usepackage{blkarray}        
\usepackage{footnote}
\makesavenoteenv{tabular}
\makesavenoteenv{table}
\usepackage{stmaryrd}
\usepackage{bm}
\usepackage{makecell}

\usepackage[square,numbers,comma]{natbib}

\hypersetup{colorlinks,citecolor=blue}

\usepackage{standalone}
 
\usepackage{pgfplots}
\usepackage{pgfplotstable}
\usetikzlibrary{pgfplots.groupplots}
\usetikzlibrary{fit}
\usepackage{filecontents}
\usepackage{tikz-3dplot}

\usepgfplotslibrary{patchplots}

\usepackage{tikz}
\usepackage{tikzscale}

\usetikzlibrary{arrows,3d}

\usetikzlibrary{shapes.geometric}
\usetikzlibrary{positioning}
\usetikzlibrary{calc}
\usetikzlibrary{decorations.pathreplacing} 
\usetikzlibrary{decorations.markings}

\usetikzlibrary{spy}
\usetikzlibrary{snakes,arrows,shapes}    
\usetikzlibrary{spy}
\usetikzlibrary{backgrounds,tikzmark}  
\usetikzlibrary{decorations}
\usetikzlibrary{positioning}
\usetikzlibrary{patterns}
\usetikzlibrary{calc} 

\usepackage{booktabs}   
\usepackage{makecell} 
\usepackage{colortbl}   

\usepackage{xspace}
\usepackage{color}     
\usepackage{setspace}   
 
\usepackage{soul}
\usepackage{ulem}
\usepackage{currfile}
\usepackage{import}

\usepackage{authoraftertitle} 
\usepackage{authblk} 

\usepackage[noabbrev]{cleveref}

\usepackage{enumitem}
\usepackage{multicol} 
\usepackage{algorithm}
\usepackage{algpseudocode}
\usepackage{algorithmicx}
\usepackage[section]{placeins}

\usepackage{environ}

    \usepackage{framed}
    
    \setlist{noitemsep,topsep=2pt,parsep=0pt,partopsep=0pt}
  
    \usepackage{lmodern}

    \usepgfplotslibrary{external} 
    \usetikzlibrary{external}
   \pgfdeclarelayer{background}
   \pgfdeclarelayer{foreground}
   \pgfsetlayers{background,main,foreground}

\pgfplotsset{compat=1.8}

\newcommand{\findmax}[3]{
    \pgfplotstablesort[sort key={#2},sort cmp={float >}]{\sorted}{#1}%
    \pgfplotstablegetelem{0}{#2}\of{\sorted}%
    \let #3=\pgfplotsretval%
}

\pgfplotscreateplotcyclelist{myCycleListBlackAndWhite}
{%
  {black, mark=o},
  {black, mark=square},
  {black, mark=triangle},
  {black, mark=diamond}, 
  {black, mark=pentagon}, 
  {black, mark=*},
  {black, mark=square*},
  {black, mark=triangle*},
  {black, mark=diamond*}, 
  {black, mark=pentagon*}, 
}

\definecolor{darkgreen}{rgb}{0,0.4,0} 
\definecolor{darkbrown}{rgb}{0.5, 0.396, 0.09}

\definecolor{c1}{rgb}{0.0, 0.4196078431372549, 0.6431372549019608}
\definecolor{c2}{rgb}{1.0, 0.5019607843137255, 0.054901960784313725}
\definecolor{c3}{rgb}{0.6705882352941176, 0.6705882352941176,
0.6705882352941176} \definecolor{c}{rgb}{0.34901960784313724, 0.34901960784313724, 0.34901960784313724}
\definecolor{c4}{rgb}{0.37254901960784315, 0.6196078431372549,
0.8196078431372549} \definecolor{c}{rgb}{0.7843137254901961, 0.3215686274509804, 0.0}
\definecolor{c5}{rgb}{0.5372549019607843, 0.5372549019607843,
0.5372549019607843} \definecolor{c}{rgb}{0.6352941176470588, 0.7843137254901961, 0.9254901960784314}
\definecolor{c6}{rgb}{1.0, 0.7372549019607844, 0.4745098039215686}
\definecolor{c7}{rgb}{0.8117647058823529, 0.8117647058823529,
0.8117647058823529}

\pgfplotscreateplotcyclelist{myCycleListColor}
{%
  {blue, mark=o},
  {red, mark=square},
  {darkbrown, mark=triangle},
  {darkgreen, mark=diamond},
  {black, mark=pentagon},
  {blue, mark=*},
  {red, mark=square*},
  {darkbrown, mark=triangle*},
  {darkgreen, mark=diamond*},
  {black, mark=pentagon*},
}

\pgfplotsset{every axis/.append style= 
              {
                font=\scriptsize,
                mark size=3,
                legend style={font=\tiny, mark size=3, draw=none, fill=none},
                legend cell align=left,
                cycle list name=myCycleListColor,
                scaled y ticks = false,
				scaled x ticks = false,
				trim axis left,
				trim axis right,
				sharp plot,
				tick label style ={font=\tiny},
				label style ={font=\scriptsize},
				very thin,
				ymajorgrids=true,
				grid style=dotted,
				legend pos= north east,
				height=\figHeight, width=\figWidth,
              }
            }

\pgfplotsset{every tick label/.append style={font=\tiny}}
\pgfplotsset{every label/.append style={font=\scriptsize}}
\pgfplotsset{every axis legend/.append style={font=\tiny}}
\pgfplotsset{every axis plot/.append style={thick}}
            
\pgfplotstableset
{
  every odd row/.style={before row={\rowcolor[gray]{0.9}}},
  every head row/.style=
  {
    before row=
    {%
      \toprule
    },
    after row=\midrule
  },
  every last row/.style=
  {
    after row=\bottomrule
  }
}

\newif\ifdrawboundingbox

\pgfkeys
{ 
  /pgf/number format/.cd,
  fixed,
  set thousands separator={\,},
  1000 sep in fractionals,
}

\newcolumntype{x}[1]{>{\centering\arraybackslash}p{#1}}
\newcommand\diag[4]{%
  \multicolumn{1}{p{#2}|}{\hskip-\tabcolsep
  $\vcenter{\begin{tikzpicture}[baseline=0,anchor=south west,inner sep=#1]
  \path[use as bounding box] (0,0) rectangle (#2+2\tabcolsep,1.5+\baselineskip);
  \node[minimum width={#2+2\tabcolsep-\pgflinewidth},
        minimum  height={\baselineskip+\extrarowheight-\pgflinewidth}] (box) {};
  \draw[line cap=round] (box.north west) -- (box.south east);
  \node[anchor=south west] at (box.south west) {#3};
  \node[anchor=north east] at (box.north east) {#4};
 \end{tikzpicture}}$\hskip-\tabcolsep}}

\newcolumntype{C}[1]{>{\centering\arraybackslash}m{#1}}
\newcolumntype{R}[1]{>{\raggedright\arraybackslash}m{#1}}
\newcolumntype{L}[1]{>{\raggedleft\arraybackslash}m{#1}}

\newcommand\restr[2]{{
		\left.\kern-\nulldelimiterspace 
		#1 
		\right|_{#2} 
}}

\newcommand\compactDots{\ifmmode\ldots\else\makebox[10cm][c]{.\hfil.\hfil.}\fi}

\def\registered{{\ooalign{\hfil\raise .00ex\hbox{\tiny R}\hfil\cr\mathhexbox20D}}}

\makeatletter
\def\namedlabel#1#2{\begingroup
	\def\@currentlabel{#2}%
	\label{#1}\endgroup
}
\makeatother

\makeatletter
\renewcommand{\todo}[2][]{\tikzexternaldisable\@todo[#1]{#2}\tikzexternalenable}
\makeatother

\usepackage{etoolbox}
\AtBeginEnvironment{mysmallmatrix}{%
\setstretch{0.9}%
\setlength{\arraycolsep}{2pt}%
}%

{%
	\begin{matrix}%
}%
	{\end{matrix}}%

{%
	\begin{matrix}%
	}%
	{\end{matrix}}%
\AtBeginEnvironment{myverysmallmatrix}{%
	\footnotesize%
	\setstretch{0.55}%
	\setlength{\arraycolsep}{1.5pt}%
}%

\newcommand{\overbar}[1]{\mkern 0.1mu\overline{\mkern-0.1mu#1\mkern-0.1mu}\mkern 0.1mu}
\makeatletter
\def\bbordermatrix#1{\begingroup \m@th
	\@tempdima 4.75\p@
	\setbox\z@\vbox{%
		\def\cr{\crcr\noalign{\kern2\p@\global\let\cr\endline}}%
		\ialign{$##$\hfil\kern2\p@\kern\@tempdima&\thinspace\hfil$##$\hfil
			&&\quad\hfil$##$\hfil\crcr
			\omit\strut\hfil\crcr\noalign{\kern-\baselineskip}%
			#1\crcr\omit\strut\cr}}%
	\setbox\tw@\vbox{\unvcopy\z@\global\setbox\@ne\lastbox}%
	\setbox\tw@\hbox{\unhbox\@ne\unskip\global\setbox\@ne\lastbox}%
	\setbox\tw@\hbox{$\kern\wd\@ne\kern-\@tempdima\left[\kern-\wd\@ne
		\global\setbox\@ne\vbox{\box\@ne\kern2\p@}%
		\vcenter{\kern-\ht\@ne\unvbox\z@\kern-\baselineskip}\,\right]$}%
	\null\;\vbox{\kern\ht\@ne\box\tw@}\endgroup}
\makeatother

\title{A projected super-penalty method for the $C^1$-coupling of multi-patch isogeometric Kirchhoff plates.}

\author[1]{Luca Coradello\thanks{luca.coradello@epfl.ch, Corresponding Author}}
\author[3]{Gabriele Loli}
\author[1,2]{Annalisa Buffa}

\affil[1]{Chair of Numerical Modelling and Simulation,
 		  \'Ecole Polytechnique F\'ed\'erale de Lausanne, Lausanne, Switzerland.}

\affil[2]{Istituto di Matematica Applicata e Tecnologie Informatiche `E. Magenes' (CNR), Pavia, Italy.}
\affil[3]{Dipartimento di Matematica `F. Casorati', Universit\`a di Pavia, Pavia, Italy.}

\newcommand{\publicationDate}{\today}

\usepackage{fancyhdr}
\date{}
\fancypagestyle{plain}{%
  \fancyhf{}%
  \fancyfoot[L]{
  \vspace{-1.25cm} 
  \footnotesize{Preprint submitted to \textit{}  \hfill
  \publicationDate
  \\
  }} }
 
\begin{document}  
\normalem
\maketitle  
  
\vspace{-1.5cm} 
\hrule 
\section*{Abstract}
This work focuses on the development of a super-penalty strategy based on the $L^2$-projection of suitable coupling terms to achieve $C^1$-continuity between non-conforming multi-patch isogeometric Kirchhoff plates. In particular, the choice of penalty parameters is driven by the underlying perturbed saddle point problem from which the Lagrange multipliers are eliminated and is performed to guarantee the optimal accuracy of the method. Moreover, by construction, the method does not suffer from locking also on very coarse meshes. 
We demonstrate the applicability of the proposed coupling algorithm to Kirchhoff plates by studying several benchmark examples discretized by non-conforming meshes. In all cases, we recover the optimal rates of convergence achievable by B-splines where we achieve a substantial gain in accuracy per degree-of-freedom compared to other choices of the penalty parameters.
\vspace{0.2cm} 
\hrule

\def\Estconst{3}
 
\vspace{0.25cm}
\noindent \textit{Keywords:} isogeometric analysis, multi-patch coupling, super-penalty method, Kirchhoff plates. 
\vspace{0.25cm}
\hrule 


\section{Introduction}\label{sec:introduction}

Isogeometric analysis (IGA), firstly introduced in~\citep{Hughes2005}, is a methodology used for the numerical discretization of partial differential equations (PDEs) based on the same building blocks used in Computer Aided Design (CAD). Indeed, in IGA, the same mathematical objects, such as B-splines and non-uniform rational B-splines (NURBS)~\citep{Piegl1995}, used for the geometrical description are employed for the numerical solution of the PDE at hand. A distinguishing feature of splines is the high regularity achievable by construction, which allows the approximation of higher-order variational problems directly in their primal, for instance Kirchhoff plates~\citep{Niiranen2017,Ning2018}, Kirchhoff-Love shells \citep{Kiendl2009,Reali2015,Kiendl2015,Kiendl2016} and the Cahn-Hilliard equation~\citep{Gomez2008}. For a detailed review of the method and its recent applications, the reader is referred to~\citep{Hughes2005,Cottrell2009,Hughes2017special}, whereas its mathematical foundations can be found in~\citep{Bazilevs2006,Buffa2014}.

Although smoothness is attained naturally within a patch, geometries of engineering relevance are in general described by multiple patches, where typically the underlying spline representations are non-conforming at the common interface. Clearly, in this scenario, a direct strong coupling between patches is not straightforward to achieve. Moreover, as in the scope of this work we are interested in the Kirchhoff plate model problem, an efficient strategy to obtain $C^1$-coupling is needed since a global $C^1$-continuity is required to obtain a well-defined bilinear form for the problem at hand.
In the literature, three methods are predominantly used to achieve the latter coupling in a weak sense and they are summarized in the following.

\noindent High-order mortar methods have been studied in~\citep{Horger2019,Hirschler2019} in the context of Kirchhoff plates and Kirchhoff-Love shells, respectively, and have been extended to a general $C^n$-coupling in~\citep{Dittmann2019}. For a detailed review in the context of isogeometric analysis, we refer to the review article~\citep{Hesch2020}. However, mortar methods leads to the formulation of a saddle point problem, where the associated Lagrange multipliers constitute additional unknowns to be solved for in the global system of equations.

\noindent Nitsche method has been analyzed in~\citep{Schillinger2016} for coupling isogeometric Kirchhoff plates in the scope of immersed methods and in~\citep{Harari2015} for imposing weakly kinematic boundary conditions for fourth-order PDEs.
Although this family of method is less sensitive to the choice of parameters compared to classical penalty approaches, their formulation requires additional consistency terms which, in the Kirchhoff problem, involve the computation of derivatives of shape functions up to order three. This adds some extra steps of complexity in the implementation and increases the overall computational cost of the coupling strategy.

\noindent Finally, penalty methods are widely used in the engineering community due to their conceptual simplicity, see the seminal work~\citep{Babuska1973}. Furthermore, they can be easily and efficiently incorporated into a numerical code, where we refer to~\citep{Kiendl2010,Apostolatos2015,Duong2017,Herrema2019} for more insights and some applications in the context of isogeometric Kirchhoff-Love shells. Nonetheless, a major drawback of this approach resides in their lack of robustness with respect to the choice of penalty parameters. Typically, the choice of penalty coefficients is problem-dependent and is based on a time-consuming, heuristic process. As noted in~\citep{Herrema2019}, on one hand, if the penalty factors are chosen too small the interface constraint is satisfied only loosely. On the other hand, if the coefficients are too high, the condition number of the resulting system matrix is negatively impacted and the convergence behavior is spoiled by spurious locking phenomena.

Our contribution falls into this realm. Inspired by the \textit{super-penalty} method studied in~\citep{Babuska73}, our goal is to introduce a simple coupling procedure for the displacement and rotation fields, respectively, for non-conforming multi-patch Kirchhoff plates, which preserves the high-order optimal convergence rates achievable by B-splines while mitigating the detrimental effects related to locking. 
To alleviate the over-constraint of the solution space we perform an $L^2$-projection of the penalty terms onto a space of reduced degree defined on the slave side of the coupling interface where, motivated by the work in~\citep{Brivadis2015} for mortar methods, we select a $p/p-2$ pairing, where $p$ denotes the B-splines degree. In particular, starting from the perturbed saddle point formulation of the Kirchhoff plate model problem, we show how the corresponding Lagrange multipliers can be eliminated from the system and, more importantly, how the perturbation gives us insights into the optimal choice for the penalty coefficients.
Indeed, the proposed methodology is truly parameter-free, as the penalty factors are fully determined by the given physical constants, the geometry and its discretization, i.e. mesh size and spline degree.  
We remark that the proposed methodology is especially advantageous for moderate degrees $p=2,3$, where locking phenomena are particularly pronounced and the $L^2$-projection proves to be an effective and computationally efficient remedy.

\noindent Then, we address the ill-conditioning issues stemming from our choice of super-penalty parameters. We adapt the block preconditioner based on an inexact Schur Complement Reduction (SCR) introduced in~\citep{Liu2019,Liu2020} and we combine it with a preconditioner tailored to the isogeometric discretization of the Kirchhoff plate, where we exploit the tensor product structure of B-splines and an efficient algorithm for the solution of the arising Sylvester-like system; for a detailed derivation we refer to~\citep{Tani2016,Montardini2018,Loli2019}.

\noindent Finally, we show through several numerical benchmarks the optimal convergence properties of the presented methodology, where our approach does not suffer from locking also on very coarse meshes. This leads to a substantial improvement in the accuracy achievable per degree-of-freedom (dof).

The structure of the paper is as follows. \Cref{sec:bsplines} provides a review of the fundamental concepts related to B-splines. \Cref{sec:formulation} describes in details the derivation of the proposed methodology and motivates our choice of penalty parameters. \Cref{sec:preconditioner} presents the ideas used in the construction of the preconditioner employed in this work. In \Cref{sec:numericalExamples} the method is validated on several numerical benchmarks and it is applied to the analysis of an idealized multi-patch design of an L-bracket. Finally, some conclusions are drawn in \Cref{sec:conclusions}.
\newcommand{\picsDir}{pictures/numericalExamples/pics}
\newcommand{\graphDir}{pictures/numericalExamples/graphs}
\newcommand{\dataDir}{pictures/numericalExamples/data}

\newtheorem{remark}{Remark}
\setlength{\extrarowheight}{0.05cm}

\section{A brief introduction to B-splines} \label{sec:bsplines}
In this section, some definitions and fundamentals related to B-splines and NURBS are reviewed. We refer the reader to~\citep{Piegl1995, Cottrell2009, Hoellig_book}, and references therein, for a comprehensive review of B-splines and their role in isogeometric analysis.

\noindent Starting from two integers $p, n$, a univariate B-spline basis function $b_{i,p}$ of degree $p$ is generated starting from a non-decreasing sequence of real values referred to as knot vector, denoted in the following as $\Xi = \left\lbrace \xi_1, \ldots, \xi_{n+p+1} \right\rbrace$.
It is worth mentioning that the smoothness of the obtained B-spline basis is $C^{p-k}$ at every knot, where $k$ denotes the multiplicity of the considered knot, while it is $C^\infty$ elsewhere. In the remainder of this work, we consider only splines of maximum continuity, i.e. $C^{p-1}$.
The definition of multivariate B-splines $\mathcal{B}_{\mathbf{i},\mathbf{p}} (\bm{\eta} )$ is achieved in a straight-forward manner using the tensor product of univariate B-splines as:
\begin{equation*}
\mathcal{B}_{\mathbf{i},\mathbf{p}} (\bm{\eta} ) = \prod_{j=1}^{\widehat{d}} b_{i_{j},p_{j}}^{j} (\eta_j)\, ,
\end{equation*}
where $\widehat{d}$ denotes the dimension of the parameter space. Additionally, the multi-index $\mathbf{i}=\left\{ i_{1},...,i_{\widehat{d}}\right\} $
denotes the position in the tensor product structure and $\mathbf{p}=\left\{ p_{1},...,p_{\widehat{d}}\right\} $ indicates the vector of polynomial degrees, associated to the corresponding parametric dimension $\bm{\eta} = \eta_1,\ldots,\eta_{\widehat{d}}\,$, respectively. 

\noindent Then, let us define a domain $\Omega \in \mathbb{R}^{d}$ described by a B-spline parametrization $\mathbf{F}$ as a linear combination of multivariate B-spline basis functions and corresponding control points as follows:
\begin{equation*}
\Omega = \mathbf{F}(\widehat{\Omega}) \quad \text{with} \quad \mathbf{F} (\bm{\eta} ) =\sum_{\mathbf{i}}\mathcal{B}_{\mathbf{i},\mathbf{p}}(\bm{\eta} ) \mathbf{P}_{\mathbf{i}} \, ,
\end{equation*}
where the coefficients $\mathbf{P}_{\mathbf{i}}\in\mathbb{R}^{d}$ of the linear combination are the control points and $d$ represents the dimensionality of the physical space. Although not treated here, it is straightforward to extend the notation to NURBS, for details see~\citep{Cottrell2009}.
In the rest of the paper, without loss of generality, the degree vector $ \mathbf{p} $ will be considered equal in each parametric direction and therefore simplified to a single scalar value $p$. Further, the vectors $\mathbf{i}$ and $\bm{\eta}$ will be omitted to simplify the notation.

\noindent Finally, we can introduce the following discrete space formed by multivariate B-splines of degree $p$:
\begin{align*}
S_{h}^p(\Omega) = \text{span }\left\lbrace b \circ \mathbf{F}^{-1} \, \vert \, b \in \mathcal{B} \right\rbrace \, .
\end{align*}

\section{The projected super-penalty method} \label{sec:formulation}

In this section, we introduce a method which alleviates locking phenomena arising when coupling non-conforming isogeometric patches. Inspired by the work presented in~\citep{Brivadis2015} in the context of isogeometric mortar methods, the proposed technique is based on the projection of the coupling terms at the interface, typically defined in terms of the degree $p$ of the solution space, onto a reduced space of B-splines of degree $p^{\text{red}} = p - 2$ defined on the slave side of the interface.  

\subsection{The strong form of the Kirchhoff plate problem}

Let us introduce the governing PDE, characterized by the bilaplace differential operator, that describes the bending-dominated problem of a Kirchhoff plate, following the notation in~\citep{Reali2015}. Let us define an open set $\Omega \subset \mathbb{R}^2$ with a sufficiently smooth boundary $\partial \Omega$, such that the normal vector $\boldsymbol{n}$ to the boundary is well-defined (almost) everywhere.   
Let us also introduce two admissible splittings of the boundary $\Gamma = \partial \Omega$ into $\Gamma =\overbar{\Gamma_u \cup \Gamma_Q}$ and $\Gamma =\overbar{\Gamma_{\phi} \cup \Gamma_{M}}$, such that $\Gamma_u \cap \Gamma_Q = \varnothing$ and $\Gamma_{\phi} \cap \Gamma_{M} = \varnothing$, respectively. Consequently, the strong form of the problem reads: 
\begin{alignat}{2}\label{eq:strongBVPPlate}
	D \varDelta^2 u &= g \quad &&\text{in} \quad \Omega \nonumber \\ 
	u &= u_{\space \Gamma} \quad &&\text{on} \quad \Gamma_u \nonumber \\ 
	- \nabla u \cdot \boldsymbol{n} &= \phi_{\space \Gamma} \quad &&\text{on} \quad \Gamma_{\phi} \nonumber \quad \\
	\nu D \varDelta u + (1 - \nu) D \, \boldsymbol{n} \cdot (\nabla \nabla u)\boldsymbol{n} &= M_{\space \Gamma} \quad &&\text{on} \quad \Gamma_{M} \nonumber \\
	D ( \nabla(\varDelta u) + (1 - \nu) \, \boldsymbol{\varPsi}(u)\,) \cdot \boldsymbol{n} &= Q_{\space \Gamma} \quad &&\text{on} \quad \Gamma_{Q}   \, ,
\end{alignat} 
where $u$ represents the deflection of the plate, $D$ its bending stiffness, $\nu$ is the Poisson ratio, $g$ is the load per unit area in the thickness direction, $u_{\space \Gamma}$, $\phi_{\space \Gamma}$, $M_{\space \Gamma}$ and $Q_{\space \Gamma}$ are the prescribed deflection, rotation, bending moments and effective shear, respectively. The bending stiffness $D$ is defined as:
\begin{align*}
D = \frac{E t^3}{12(1 - \nu^2)} \, ,
\end{align*}
where $E$ is the Young modulus and $t$ denotes the thickness of the plate. For the sake of simplicity and without loss of generality, these are assumed to be a constant in $\Omega$.
Finally, the differential operator $\boldsymbol{\varPsi}(\cdot)$ reads:
\begin{align*}
\boldsymbol{\varPsi}(\cdot) = \left[
\frac{\partial^3 (\cdot)}{\partial x \partial^2 y} \, , \, \frac{\partial^3 (\cdot)}{\partial^2 x \partial y} \right]^\top \, .
\end{align*}

\subsection{The multi-patch formulation of the perturbed saddle point Kirchhoff problem}

\begin{figure}
	\centering
	\includegraphics[width=0.65\textwidth]{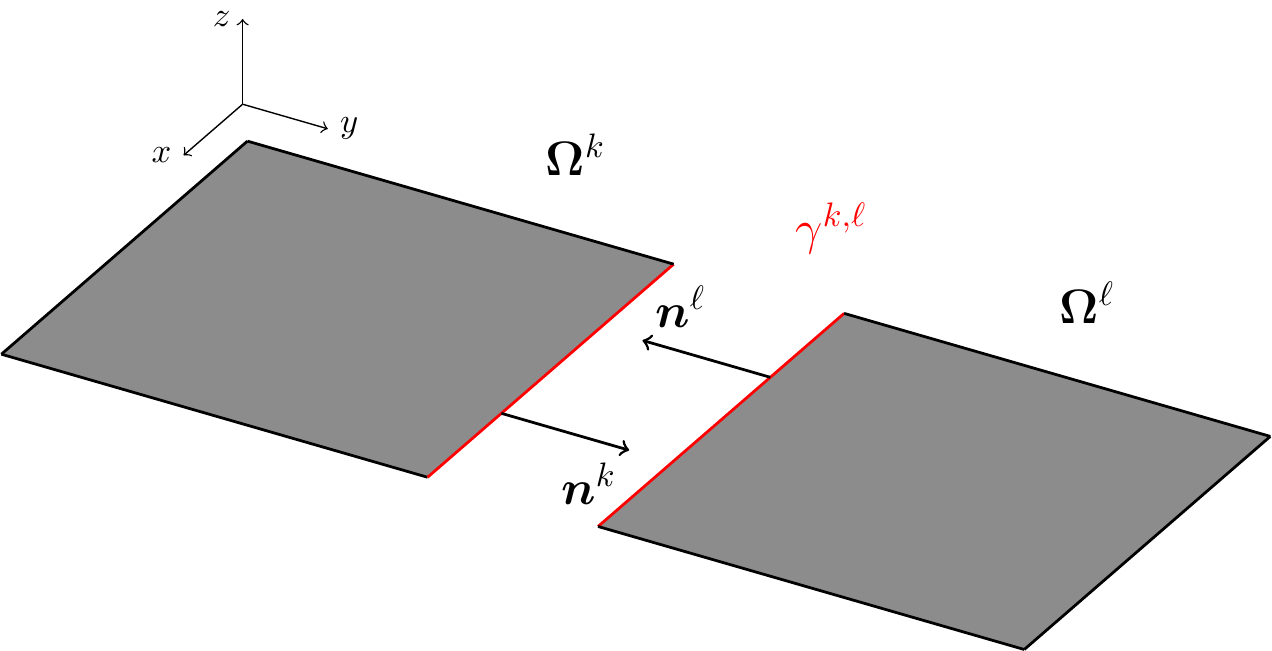}
	\caption{Example of two subdomains $\Omega^k, \Omega^\ell$ with their coupling interface $\gamma^{k,\ell}$, highlighted in red, and their corresponding normal vectors $\bm{n}^k, \bm{n}^\ell$. Note that we have separated the subdomains for visualization purposes. For a correct interpretation of the colors, the reader is referred to the web version of this manuscript.}\label{fig:coupling_plate}
\end{figure}

Here, following the notation used in~\citep{Brivadis2015}, we introduce a decomposition of $\Omega$ into $N$ non-overlapping subdomains $\Omega^i$ such that:
\begin{align*}
\overbar{\Omega} = \bigcup_{i=1}^N \overbar{\Omega^i} \, , \quad \text{where} \quad \Omega^i \cap \Omega^j = \varnothing \quad \text{for} \quad i \neq j \, .
\end{align*}
Now, let us define the interface $\gamma^{k,\ell}$ between two adjacent patches $\Omega^{k}, \Omega^{\ell}, 1 \leq k, \ell \leq N$ as the intersection of their corresponding boundaries:
\begin{align*}
\gamma^{k,\ell} = \partial \Omega^{k} \cup \partial \Omega^{\ell} \, .
\end{align*}
Then, the skeleton $\Gamma$ is defined as the union of all non-empty interfaces (which we suppose to be labeled with an index $\ell=1, \ldots, L$) and reads:
\begin{align*}
\Gamma = \bigcup_{\ell=1}^L \gamma^\ell \, .
\end{align*} 
Consequently, we can denote by $u^k$ and $\bm{n}^k$ the value of the primary field and the outward normal on $\partial \Omega^k$, and $u^\ell$ and $\bm{n}^\ell$ the value of the primary field and outward normal on the neighboring subdomain $\partial \Omega^\ell$, see~\Cref{fig:coupling_plate} for an example on two patches.
Then, for each interface $\gamma^{k,\ell}$ we can write the following coupling conditions:
\begin{align*}
u^k - u^\ell &= 0 \quad \text{on} \quad \gamma^{k,\ell}  \\
\nabla u^k \cdot \boldsymbol{n}^k + \nabla u^\ell \cdot \boldsymbol{n}^\ell &= 0 \quad \text{on} \quad \gamma^{k,\ell} \, ,
\end{align*}
which can be rewritten using the standard jump and normal jump operators, respectively, as:
\begin{align}\label{eq:coupling_conds}
\llbracket u \rrbracket &= 0 \quad \text{on} \quad \gamma^{k,\ell} \nonumber \\
\llbracket \nabla u \rrbracket _n &= 0 \quad \text{on} \quad \gamma^{k,\ell} \, . 
\end{align}
Further, given $1 \leq s, t \leq L$, $s \neq t$, we denote the cross-points by $c^{s,t}=\overline{\gamma}^s \cap \overline{\gamma}^t$ and we label them with an ordered index $c^s$, $s=1,\ldots,S$.
For ease of notation and without loss of generality, in the following we assume the flexural rigidity $D$ to be constant in $\Omega$ and the Poisson ratio $\nu$ to be zero. Further, we assume that the values prescribed as natural boundary conditions are zero as well. 

Now, let us introduce for each subdomain $\Omega^i$ the following space:
\begin{align*}
H^2_\star(\Omega^i)= \left\lbrace v^i \in H^2(\Omega^i): \, v^i |_{\partial\Omega \cap \partial \Omega^i} = \frac{\partial v^i}{\partial n} \bigg\rvert_{\partial\Omega \cap \partial \Omega^i} =0  \right\rbrace \, ,
\end{align*}
from which the following broken Sobolev space can be characterized as:
\begin{equation*}
	V=\left\lbrace v \in L^2(\Omega) :\,   v |_{\Omega^i} \in H^2_\star(\Omega^i),\,i=1,\ldots, N \, , \: v \text{ is continuous in } c^s, \:s=1,\ldots,S	\right\rbrace \, ,
\end{equation*}
endowed with the broken norm $\| \cdot \|^2_{V}= \sum_{i=1}^{N} \| \cdot \|^2_{H^2(\Omega^i)}$.
Then, let us also define the spaces:
\begin{align*}
	H^{\frac{1}{2}}_\star(\Gamma)&=\left\lbrace \left\llbracket \nabla v \right\rrbracket_n :\, v \in V \right\rbrace  \\
	H^{\frac{3}{2}}_\star(\Gamma)&=\left\lbrace \llbracket v \rrbracket:\, v \in V \right\rbrace \, .
\end{align*}
Lastly, we need to introduce the following dual spaces:
\begin{align*}
 	&Q_1=\left[H^{\frac{3}{2}}_\star(\Gamma)\right]'  \\
 	&Q_2=\left[H^{\frac{1}{2}}_\star(\Gamma)\right]'\, . 
\end{align*}
We are now ready to formulate~\eqref{eq:strongBVPPlate} as a saddle point problem. Given $f \in V'$, find $(u,\lambda_1,\lambda_2) \in V \times Q_1 \times Q_2$ such that:
\begin{align}\label{eq:continuous_mortar_problem}
\sum_{i=1}^{N}\int_{\Omega^i} D \, \nabla (\nabla u ) : \nabla (\nabla v ) + \sum_{\ell=1}^{L}\left(\int_{\gamma^{\ell}} \llbracket v \rrbracket \lambda_1 + \int_{\gamma^{\ell}} \left\llbracket \nabla v \right\rrbracket_n \lambda_2 \right) &= (f,v) &\qquad  &\forall v \in V \nonumber \\
\sum_{\ell=1}^{L}\int_{\gamma^{\ell}} \llbracket u \rrbracket \mu_1 &=0 &\qquad &\forall \mu_1 \in Q_1 \nonumber \\
\sum_{\ell=1}^{L}\int_{\gamma^{\ell}} \left\llbracket \nabla u \right\rrbracket_n \mu_2 &=0 &\qquad &\forall \mu_2 \in Q_2  \, .
\end{align}
We also define three continuous bilinear forms $a:V \times V \rightarrow \mathbb{R}$, $b_1:V \times Q_1 \rightarrow \mathbb{R}$ and $b_2:V \times Q_2 \rightarrow \mathbb{R}$ as follows:
\begin{align*}
a(u,v) &= \sum_{i=1}^{N}\int_{\Omega^i} D \, \nabla (\nabla u ) : \nabla (\nabla v ) &\qquad &u,v \in V \\
b_1(v,\mu_1) &= \sum_{\ell=1}^{L}\int_{\gamma^{\ell}} \llbracket v \rrbracket \mu_1 &\qquad &v \in V, \, \mu_1 \in Q_1 \\
b_2(v,\mu_2) &= \sum_{\ell=1}^{L}\int_{\gamma^{\ell}} \left\llbracket \nabla v \right\rrbracket_n \mu_2 &\qquad &v \in V, \, \mu_2 \in Q_2 \, .
\end{align*}
Now, given $\varepsilon^{(\ell)}_1, \varepsilon^{(\ell)}_2 > 0, \ell = 1,\ldots, L$, we can introduce the singularly perturbed version of~\eqref{eq:continuous_mortar_problem}:
given $f \in V'$, find $(u_{\varepsilon},\lambda_{1,\varepsilon},\lambda_{2,\varepsilon}) \in V \times L^2(\Gamma) \times L^2(\Gamma)$, such that
\begin{align}\label{eq:continuous_perturbed_problem}
\sum_{i=1}^{N} \int_{\Omega^i} D \, \nabla (\nabla u_{\varepsilon} ) : \nabla (\nabla v ) + \sum_{\ell=1}^{L}\left(\int_{\gamma^{\ell}} \llbracket v \rrbracket \lambda_{1,\varepsilon} + \int_{\gamma^{\ell}} \left\llbracket \nabla v \right\rrbracket_n \lambda_{2,\varepsilon}\right) &= (f,v) & \qquad &\forall v \in V \nonumber \\
\sum_{\ell=1}^{L} \left( \int_{\gamma^{\ell}} \llbracket u_{\varepsilon} \rrbracket \mu_1 - \varepsilon^{(\ell)}_1 \int_{\gamma^{\ell}} \lambda_{1,\varepsilon} \mu_1 \right) &= 0  & \qquad &\forall \mu_1 \in L^2(\Gamma) \nonumber \\
\sum_{\ell=1}^{L} \left( \int_{\gamma^{\ell}} \left\llbracket \nabla u_{\varepsilon} \right\rrbracket_n \mu_2- \varepsilon^{(\ell)}_2 \int_{\gamma^{\ell}} \lambda_{2,\varepsilon} \mu_2 \right) &=0 & \qquad &\forall \mu_2 \in L^2(\Gamma) \, . 
\end{align}
Under suitable regularity assumptions, we can provide an estimation of the error introduced by the perturbations $\varepsilon^{(\ell)}_1$ and $\varepsilon^{(\ell)}_2$ on the solution of the original saddle point problem~\eqref{eq:continuous_mortar_problem} as~\citep[Remark 4.13.14]{Boffi2013}:
\begin{align}\label{eq:perturbation_est}
\vert \vert u - u_{\varepsilon} \vert \vert_V + \vert \vert \lambda_1 - \lambda_{1,\varepsilon} \vert \vert_{L^2(\Gamma)} + \vert \vert \lambda_2 - \lambda_{2,\varepsilon} \vert \vert_{L^2(\Gamma)} &\leq C \left[ \tilde{\varepsilon}_{1}\left(\sum_{\ell=1}^{L}  \|\lambda_1\|^2_{H^{3/2}(\gamma^{\ell})} \right)^{\frac{1}{2}} \right. \nonumber \\ 
&+ \left. \tilde{\varepsilon}_{2} \left(\sum_{\ell=1}^{L}   \|\lambda_2\|^2_{H^{1/2}(\gamma^{\ell})} \right)^{\frac{1}{2}}\right] \, 
\end{align}
where we have defined:
\begin{align*}
\tilde{\varepsilon}_{1}=\max_{\ell=1,\ldots,L} \varepsilon^{(\ell)}_1, \quad \text{ and }\quad \tilde{\varepsilon}_{2}=\max_{\ell=1,\ldots,L} \varepsilon^{(\ell)}_2 \, .
\end{align*}

\subsection{The projected super-penalty formulation}
\begin{figure}
	\centering
	\begin{subfigure}[t]{\textwidth}
	\centering
		\includegraphics[width=0.6\textwidth]{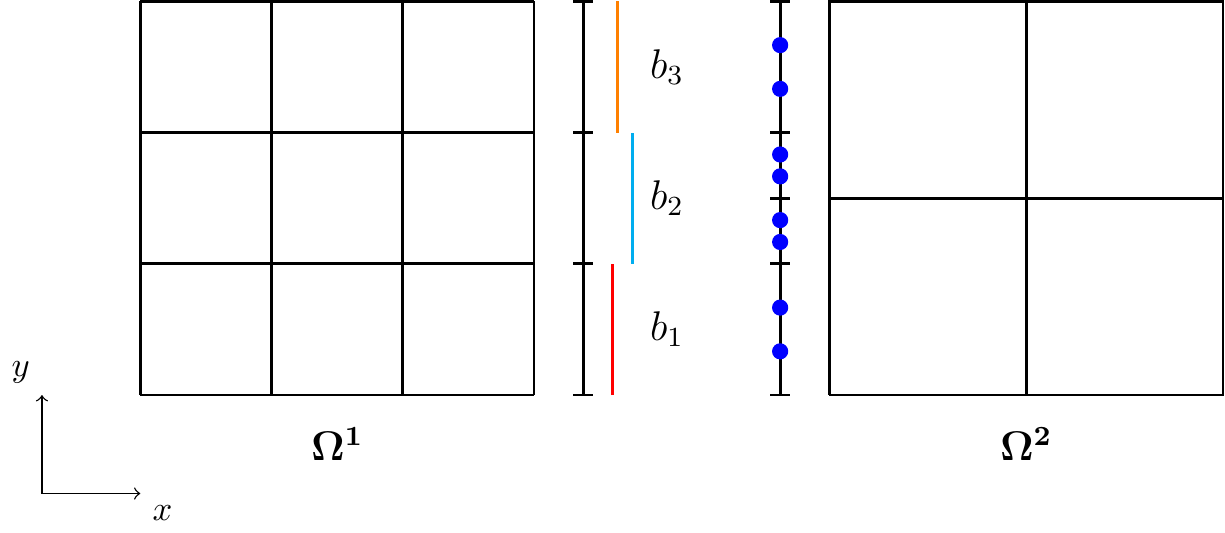}
	\caption{Interface splines of degree $p^{\text{red}} = p - 2 = 0$ associated to a reduced knot vector on the slave side \mbox{$\Xi_\star^{\ell} = \left[ 0 \,\, 1/3 \,\, 2/3 \,\, 1 \right]$} and corresponding intersection mesh for integration.}
	\end{subfigure}
	\hfill
	\begin{subfigure}[t]{\textwidth}
	\centering
		\includegraphics[width=0.6\textwidth]{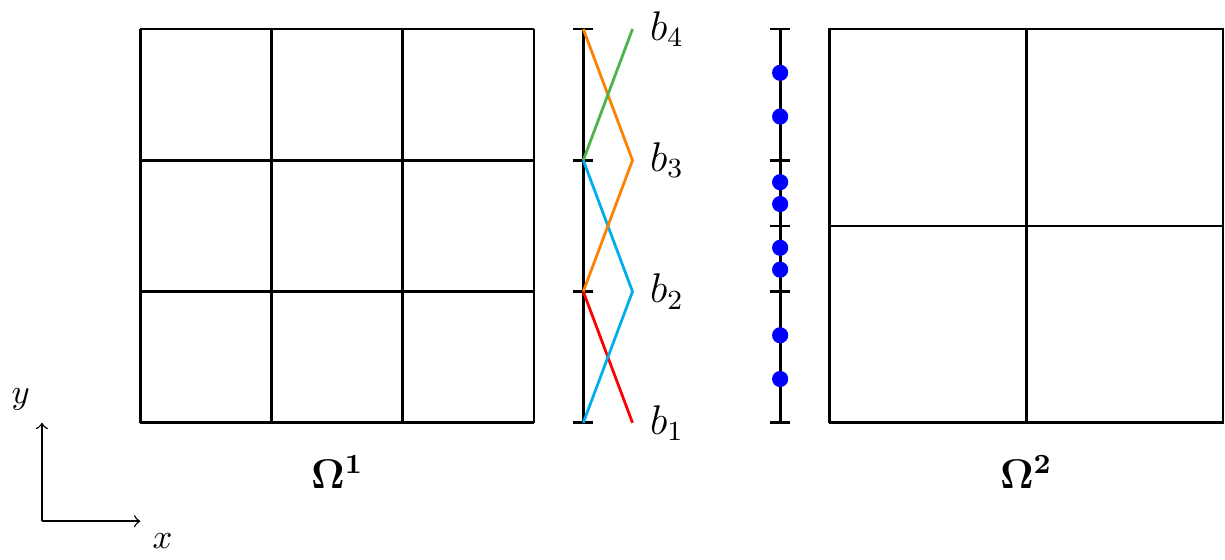}
	\caption{Interface splines of degree $p^{\text{red}} = p - 2 = 1$ associated to a reduced knot vector on the slave side \mbox{$\Xi_\star^{\ell} = \left[ 0 \, 0 \,\, 1/3 \,\, 2/3 \,\, 1 \, 1 \right]$} and corresponding intersection mesh for integration.}
	\end{subfigure}
	\caption{Example of the projection setup on a coupling interface. We select the finer side (on $\Omega^1$ in this example) to define the reduced space for the projection. Additionally, an intersection mesh at the interface is created only for integration purposes to properly compute the projected penalty terms.}\label{fig:proj_setup}
\end{figure}

For each patch $\Omega^i$, we assume $p \geq 2$ and we indicate with $\overline{\mathcal{S}^p_{h}(\Omega^i)}$ the space trivially obtained extending by zero the elements of $\mathcal{S}^p_{h}(\Omega^i)$ over $\Omega \setminus \Omega^i$. Additionally, let us define:
\begin{align*}
X_{i,h}=\mathrm{span} \{ b \in \overline{\mathcal{S}^p_{h}(\Omega^i)} \}.
\end{align*} 
Consequently, let us denote by $V_{i,h} \subset X_{i,h}$ the finite-dimensional space given by the span of B-splines defined on the corresponding subdomain $\Omega^i$, where the exact characterization of $V_{i,h}$ depends on the chosen boundary conditions, for further details we refer to~\citep{Ciarlet2002}.
This allows us to introduce the following finite dimensional subspace of $V$,
\begin{align*}
V_h &= \left\lbrace v \in \bigcup_{i=1}^{N} V_{i,h}:\: v \text{ is continuous in } c^s, \:s=1,\ldots,S \right\rbrace \, .
\end{align*}
Moreover, for each interface $\gamma^{\ell}$, we denote by $\Xi^{\ell}$ the knot vector on $\gamma^{\ell}$ inherited from the slave side. 
Motivated by the choice of the $p/p-2$ stable pairing in~\citep{Brivadis2015}, we construct the following isogeometric space $\mathcal{S}^{p-2}_{h}(\gamma^{\ell})$ on the reduced knot vector $\Xi_\star^{\ell}$ obtained by removing from $\Xi^{\ell}$ the first and last two knots, where an example is depicted in~\Cref{fig:proj_setup} for $p=2,3$. Similarly to before, we indicate with $\overline{\mathcal{S}^{p-2}_{h}(\gamma^{\ell})}$ the space obtained extending by zero over $\Gamma \setminus \gamma^{\ell}$ the elements of ${\mathcal{S}^{p-2}_{h}(\gamma^{\ell})}$.
We can now define the discrete counterpart of the Lagrange multiplier spaces as:
\begin{align*}
	Q_h = Q_{1,h} &= Q_{2,h} = \bigcup_{\ell=1}^{L} \overline{\mathcal{S}^{p-2}_{h}(\gamma^{\ell})} \, .
\end{align*} 
With these definitions at hand, the discretized version of~\eqref{eq:continuous_perturbed_problem} reads: find $\left( u_h, \lambda_{1,h},\lambda_{2,h}  \right) \in V_h \times Q_h \times Q_h$ such that:
\begin{align}\label{eq:discrete_perturbed_problem}
\sum_{i=1}^{N} \int_{\Omega^i} D \, \nabla (\nabla u_{h} ) : \nabla (\nabla v_h ) + \sum_{\ell=1}^{L}\left(\int_{\gamma^{\ell}} \llbracket v_h \rrbracket \lambda_{1,h} + \int_{\gamma^{\ell}} \left\llbracket \nabla v_h \right\rrbracket_n \lambda_{2,h}\right) &= (f,v_h) & \qquad &\forall v_h \in V_h  \nonumber \\
\sum_{\ell=1}^{L} \left( \int_{\gamma^{\ell}} \llbracket u_{h} \rrbracket \mu_{1,h} - \frac{1}{\alpha_{\text{defl}}^\ell} \int_{\gamma^{\ell}} \lambda_{1,h} \mu_{1,h} \right) &= 0  & \qquad &\forall \mu_{1,h} \in Q_h \nonumber \\
\sum_{\ell=1}^{L} \left( \int_{\gamma^{\ell}} \left\llbracket \nabla u_{h} \right\rrbracket_n \mu_{2,h} - \frac{1}{\alpha_{\text{rot}}^\ell} \int_{\gamma^{\ell}} \lambda_{2,h} \mu_{2,h} \right) &=0 & \qquad &\forall \mu_{2,h} \in Q_h \, , 
\end{align}
where $\alpha_{\text{defl}}^\ell$ and $\alpha_{\text{rot}}^\ell$ are ``large'' parameters associated to the deflections and rotations, respectively. In general, they depend on the problem definition, e.g. the physical constant $D$, the mesh size and spline degree, where a full characterization of our choice will be given later in the section.
We can now formally eliminate the Lagrange multipliers and recast~\eqref{eq:discrete_perturbed_problem} into its primal form. Indeed, we can write:
\begin{align*}
{\lambda_{1,h}}|_{\gamma^{\ell}} &= \alpha_{\text{defl}}^\ell \Pi^{\ell} \llbracket u_h \rrbracket  \\ 
{\lambda_{2,h}}|_{\gamma^{\ell}} &= \alpha_{\text{rot}}^\ell \Pi^{\ell} \left\llbracket \nabla u_h \right\rrbracket_n \, , 
\end{align*}
where $\Pi^{\ell}:L^2(\gamma^{\ell}) \rightarrow \mathcal{S}^{p-2}_{h}(\gamma^{\ell})$ denotes the $L^2$-projection, associated to the interface $\gamma^{\ell}$, onto the reduced space $\mathcal{S}^{p-2}_{h}(\gamma^{\ell})$.
Finally, employing the previous results and the properties of the $L^2$-projection, the resulting discretized bilinear form, augmented by suitable penalty terms that weakly enforce the coupling conditions~\eqref{eq:coupling_conds}, reads: find $u_h  \in \: V_h$ such that
\begin{align}\label{eq:reduced_problem}
&\sum_{i=1}^{N}\int_{\Omega^i} D \, \nabla (\nabla u_h ) : \nabla (\nabla v_h ) + \nonumber \\
&+ \sum_{l=1}^{L} \left( \alpha_{\text{defl}}^\ell \int_{\gamma^{\ell}} \Pi^{\ell} \llbracket u_h \rrbracket \Pi^{\ell} \llbracket v_h \rrbracket + \alpha_{\text{rot}}^\ell \int_{\gamma^{\ell}} \Pi^{\ell} \left\llbracket \nabla u_h \right\rrbracket_n \Pi^{\ell} \left\llbracket \nabla v_h \right\rrbracket_n \right) = (f,v_h) \qquad \forall v_h \in V_h \, .
\end{align}

\subsubsection{Inf-sup test}
\begin{figure}
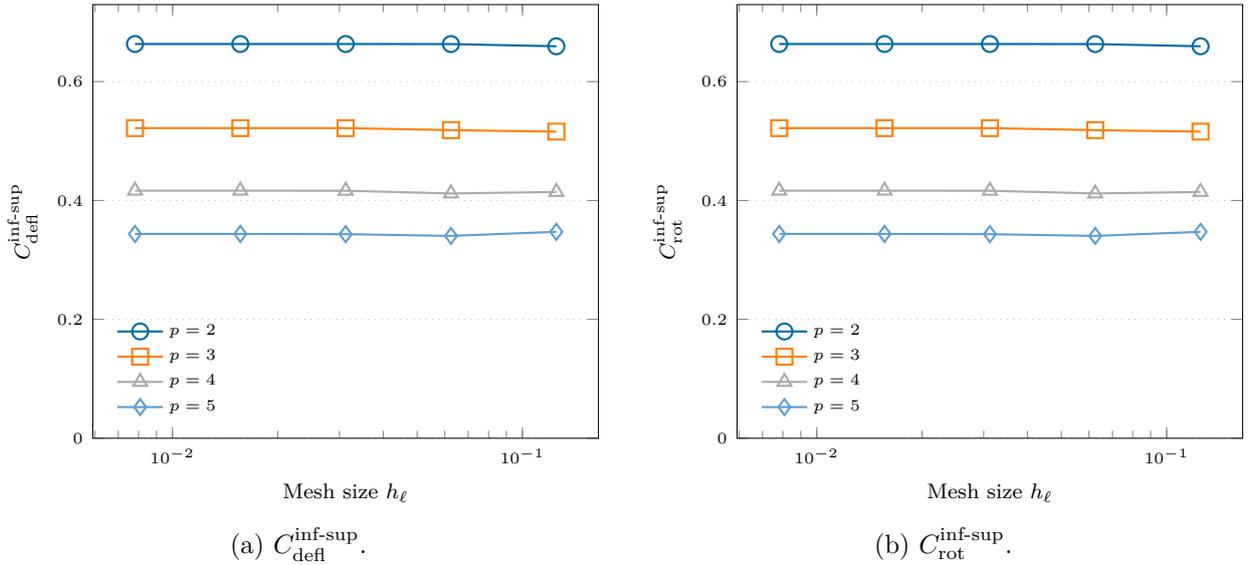

	\centering
	\begin{subfigure}[t]{0.495\textwidth}
		\centering
	\input{\graphDir/convergence_inf_sup_defl.tex}
		\caption{$C_{\text{defl}}^{\text{inf-sup}}$.}
	\end{subfigure}
	\hfill
	\begin{subfigure}[t]{0.495\textwidth}
		\centering
	\input{\graphDir/convergence_inf_sup_rot.tex}
		\caption{$C_{\text{rot}}^{\text{inf-sup}}$.}
	\end{subfigure}	
	\caption{Convergence plot of the numerical inf-sup test for $C_{\text{defl}}^{\text{inf-sup}}$ and $C_{\text{rot}}^{\text{inf-sup}}$, respectively, on different uniformly refined meshes $h_\ell = 1/2^k , k=3,\ldots,7$ and spline degrees $p=2,\ldots,5$.}
\label{fig:inf_sup}
\end{figure}

The well-posedness of~\eqref{eq:discrete_perturbed_problem}, independently of the value of the parameters $\alpha_{\text{defl}}^\ell$ and $\alpha_{\text{rot}}^\ell$, relies on the well-posedness of the underlying unperturbed problem, i.e. the problem corresponding to~\eqref{eq:discrete_perturbed_problem} where we set $\alpha_{\text{defl}}^\ell = \alpha_{\text{rot}}^\ell = +\infty \, , \ell = 1,\ldots,L$.
Although a rigorous proof of the inf-sup stability of such unperturbed problem is currently under investigation~\citep{Buffa2020}, we assess the behavior of the numerical inf-sup test for a domain $\Omega$ subdivided along a straight interface into two subdomains $\Omega^i, i=1,2$. As we are dealing with a double saddle point problem, we compute two different inf-sup constants $C_{\text{defl}}^{\text{inf-sup}}$ and $C_{\text{rot}}^{\text{inf-sup}}$, corresponding to the deflection and rotation jumps, respectively.
In the following we report the results for different discretization sizes of the interface $h_\ell = 1/2^k , k=3,\ldots,7$ and B-spline degrees $p=2,\ldots,5$, where $h_\ell$ denotes the maximum mesh size associated to the interface $\gamma^\ell$.
The numerical values of $C_{\text{defl}}^{\text{inf-sup}}$ and $C_{\text{rot}}^{\text{inf-sup}}$ are summarized in~\Cref{table:inf_sup} and are depicted for clarity in~\Cref{fig:inf_sup}. In all cases we observe that the inf-sup constants converge to some values bounded away from zero, numerically suggesting that the method is inf-sup stable.

%

\begin{table}[h]
\begin{subtable}{\textwidth}
\begin{center}
    \begin{tabular}{| x{0.5cm} | c | c | c | c | c |}
    \hline
    \diag{.1em}{.5cm}{\scalebox{0.8}{$p$}}{\scalebox{0.8}{$h_\ell$}} & $1/2^3$ & $1/2^4$ &  $1/2^5$ & $1/2^6$  & $1/2^7$ \\ \hline
    2 & 0.6594 & 0.6631 & 0.6632 & 0.6632 & 0.6632 \\ \hline    
    3 & 0.5158 & 0.5183 & 0.5217 & 0.5217 & 0.5217 \\ \hline    
    4 & 0.4144 & 0.4118 & 0.4164 & 0.4165 & 0.4165 \\ \hline    
    5 & 0.3473 & 0.3404 & 0.3433 & 0.3437 & 0.3437 \\ \hline    
    \hline
    \end{tabular}
    \caption{$C_{\text{defl}}^{\text{inf-sup}}$.}\label{table:inf_sup1}
\end{center}
\end{subtable}
\begin{subtable}{\textwidth}
\begin{center}
    \begin{tabular}{| x{0.5cm} | c | c | c | c | c |}
    \hline
    \diag{.1em}{.5cm}{\scalebox{0.8}{$p$}}{\scalebox{0.8}{$h_\ell$}} & $1/2^3$ & $1/2^4$ &  $1/2^5$ & $1/2^6$  & $1/2^7$ \\ \hline
    2 & 0.6594 & 0.6631 & 0.6632 & 0.6632 & 0.6632 \\ \hline    
    3 & 0.5158 & 0.5183 & 0.5217 & 0.5217 & 0.5217 \\ \hline    
    4 & 0.4144 & 0.4118 & 0.4164 & 0.4165 & 0.4165 \\ \hline    
    5 & 0.3473 & 0.3404 & 0.3433 & 0.3437 & 0.3437 \\ \hline    
    \hline
    \end{tabular}
    \caption{$C_{\text{rot}}^{\text{inf-sup}}$.}\label{table:inf_sup2}
\end{center}
\end{subtable}
    \caption{Results of the numerical inf-sup test for $C_{\text{defl}}^{\text{inf-sup}}$ and $C_{\text{rot}}^{\text{inf-sup}}$, respectively, on different uniformly refined meshes \mbox{$h_\ell = 1/2^k , k=3,\ldots,7$} and spline degrees $p=2,\ldots,5$.}\label{table:inf_sup}
\end{table}

\subsubsection{Coercivity test}
Then, we also assess numerically the behavior of the coercivity constant on an example with four patches $\Omega^i, i=1,\ldots,4$ separated by four straight interfaces meeting at a cross-point. In particular, we want to compute the biggest $\alpha_0$ such that:
\begin{align*}
\alpha_0 \vert \vert v_0 \vert \vert_{H^2(\Omega)} \leq a(v_0, v_0) \quad \forall v_0 \in K = \text{ker}(B_1) \cap \text{ker}(B_2) \, ,
\end{align*}
where $B_1$ and $B_2$ are the linear operators associated to the bilinear forms $b_1$ and $b_2$, respectively.
The results for different discretization sizes of the interface $h_\ell = 1/2^k , k=2,\ldots,5$ and B-spline degrees $p=2,\ldots,4$ are presented in~\Cref{table:coercivity}, from which we can numerically infer that the method is coercive on the intersection kernel.
\begin{table}[h]
\begin{center}
    \begin{tabular}{| x{0.5cm} | c | c | c | c |}
    \hline
    \diag{.1em}{.5cm}{\scalebox{0.8}{$p$}}{\scalebox{0.8}{$h_\ell$}}  & $1/2^2$ &  $1/2^3$ & $1/2^4$  & $1/2^5$ \\ \hline
    2 & 0.8049 & 0.8043 & 0.8041 & 0.8041 \\ \hline    
    3 & 0.8040 & 0.8040 & 0.8040 & 0.8040 \\ \hline    
    4 & 0.8040 & 0.8040 & 0.8040 & 0.8040 \\ \hline    
    \hline
    \end{tabular}
    \caption{Numerical estimation of the coercivity constant $\alpha_0$ on different uniformly refined meshes \mbox{$h_\ell = 1/2^k , k=2,\ldots,5$} and spline degrees $p=2,3,4$.}\label{table:coercivity}
\end{center}
\end{table}

\begin{remark}
The inf-sup and coercivity tests are performed on a reduced version of the knot vector $\Xi_\star^{\ell}$, where also the first and last internal knots of $\Xi^{\ell}$ are eliminated. This is justified by our preliminary mathematical analysis, where this choice is required. However, from a numerical standpoint, we retain the optimality of the method without performing such a reduction and in all our examples we directly employ $\Xi_\star^{\ell}$ to define the projection spaces.
\end{remark}

\subsubsection{On the choice of penalty parameters}

It is well-known that the penalized problem~\eqref{eq:reduced_problem} is variationally consistent only in the limit $\alpha_{\text{defl}}^\ell = \alpha_{\text{rot}}^\ell \rightarrow \infty \, \, \ell = 1,\ldots,L$. 
On the other hand, the well-posedness of this problem is robust with respect to the choice of the parameters $\alpha_{\text{defl}}^\ell$ and $\alpha_{\text{rot}}^\ell$. Therefore, the proposed methodology will not suffer from locking for any choice of penalty values. As a consequence, $\alpha_{\text{defl}}^\ell$ and $\alpha_{\text{rot}}^\ell$ can be chosen solely to guarantee the optimal accuracy of the method.
\begin{remark}
A clear trade-off of this choice is the negative impact on the conditioning of the resulting system matrix. A possible remedy based on an ad-hoc preconditioner will be discussed in a later section. Another drawback consists in the loss of significant digits due to the (potentially big) difference in magnitude between the penalty contribution and the internal stiffness. For this reason (amongst other which will be pointed out in the rest of manuscript), we advise to use this method in combination with splines of degree $p=2,3$, as these round-off errors occur below a tolerance threshold of significance to most engineering applications.
\end{remark}
\noindent Inspired by the method proposed in~\citep{Herrema2019} in the context of Kirchhoff-Love shells, we want to develop a fully parameter-free penalty method.
To this end, we scale the deflection and rotation penalty parameters by the physical constants, the local mesh size and the geometry as:
\begin{align}\label{eq:projectedTerms}
\alpha_{\text{defl}}^\ell &= \text{meas}( \gamma^{\ell} )^{\beta-1} \frac{E t}{(h_\ell)^\beta (1 - \nu^2)}  \nonumber\\
\alpha_{\text{rot}}^\ell &= \text{meas}( \gamma^{\ell} )^{\beta-1} \frac{E t^3}{12 (h_\ell)^\beta (1 - \nu^2)} \, ,
\end{align} 
where the exponent $\beta$ is chosen to ensure the optimal convergence of the method with respect to the degree $p$ of the underlying discretization. 
Note that all of these parameters are known and depend only on the problem definition, meaning that no user-defined factor is required.
We highlight that our choice is based on the fact that the perturbations introduced in~\eqref{eq:continuous_perturbed_problem} cannot be ``big'' compared to the accuracy with which we want to solve the original problem and the estimate provided in~\eqref{eq:perturbation_est} guides the choice of $\beta$. Moreover, as we want to recover optimal rates of convergence for the error, the exponent $\beta$ must be a function of the underlying splines degree $p$.

\noindent From the numerical experiments conducted thus far, the scaling factor $\beta = p-1$ in~\eqref{eq:projectedTerms} is necessary to ensure optimal convergence of the method in the $H^2$ norm, whereas for a scaling of $\beta = p$ we observed optimality in the $H^2$ and $H^1$ norms. Finally, a factor of $\beta = p+1$
provides optimality in the $H^2$, $H^1$ and $L^2$ norms. If not stated otherwise, we will use $\beta = p+1$ in all our numerical examples.
\begin{remark}
Although a rigorous mathematical proof of the method and the optimal choice of $\beta$ are currently under development~\citep{Buffa2020}, we believe that this allows for some extra flexibility in the proposed methodology, where the suitable scaling factor can be chosen with respect to the corresponding quantity of interest.
\end{remark}

\subsubsection{Cross-points modification}
\begin{figure}
	\centering
	\begin{subfigure}[t]{0.495\textwidth}
	\centering
		\includegraphics[width=\textwidth]{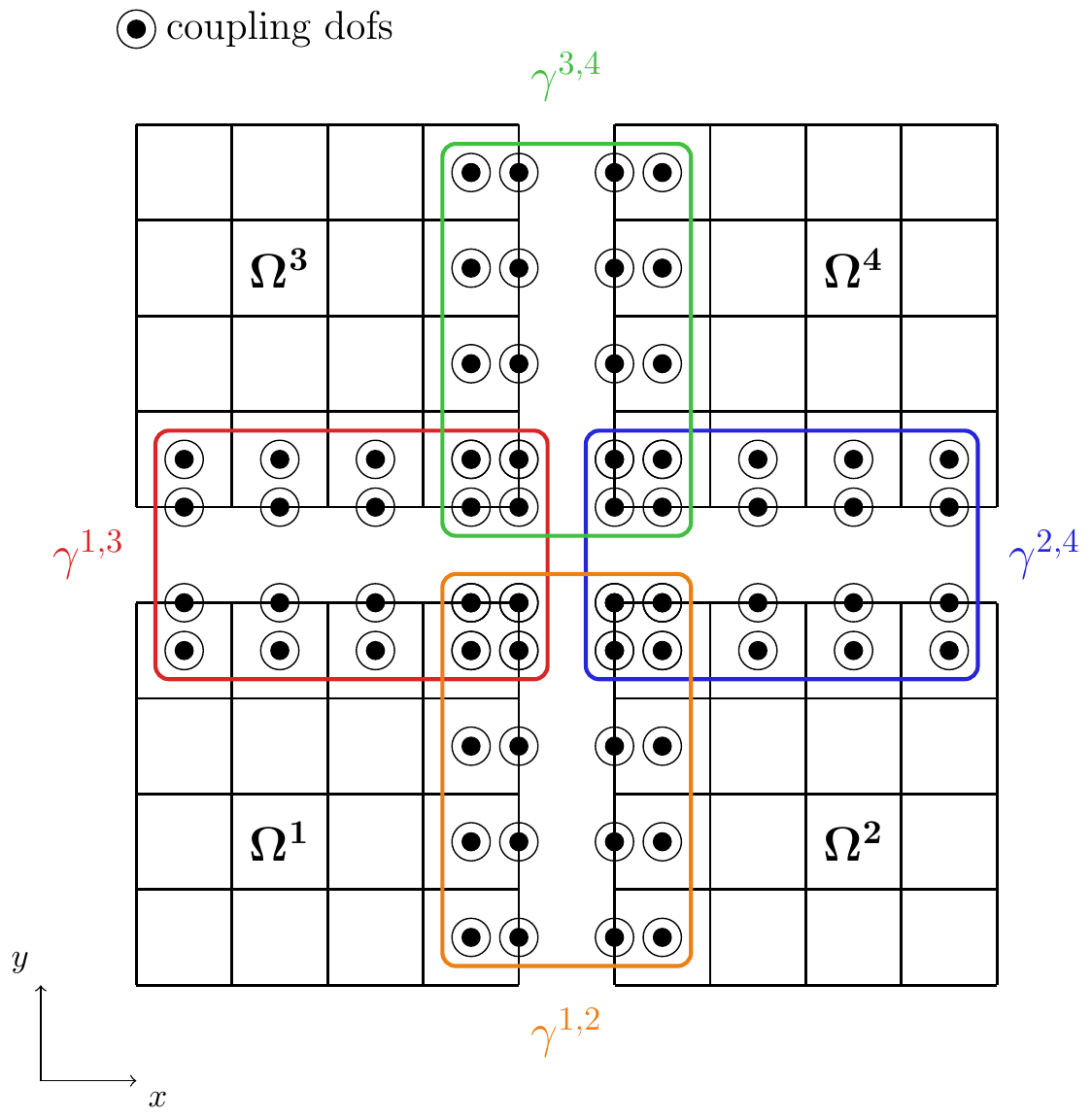}
	\caption{B-splines involved in the computation of the coupling terms, where each color colored box contains the dofs associated to the corresponding coupling interface $\gamma^{k,\ell}$.}\label{fig:dofs_interfaces}
	\end{subfigure}
	\hfill
	\begin{subfigure}[t]{0.495\textwidth}
		\centering
		\includegraphics[width=\textwidth]{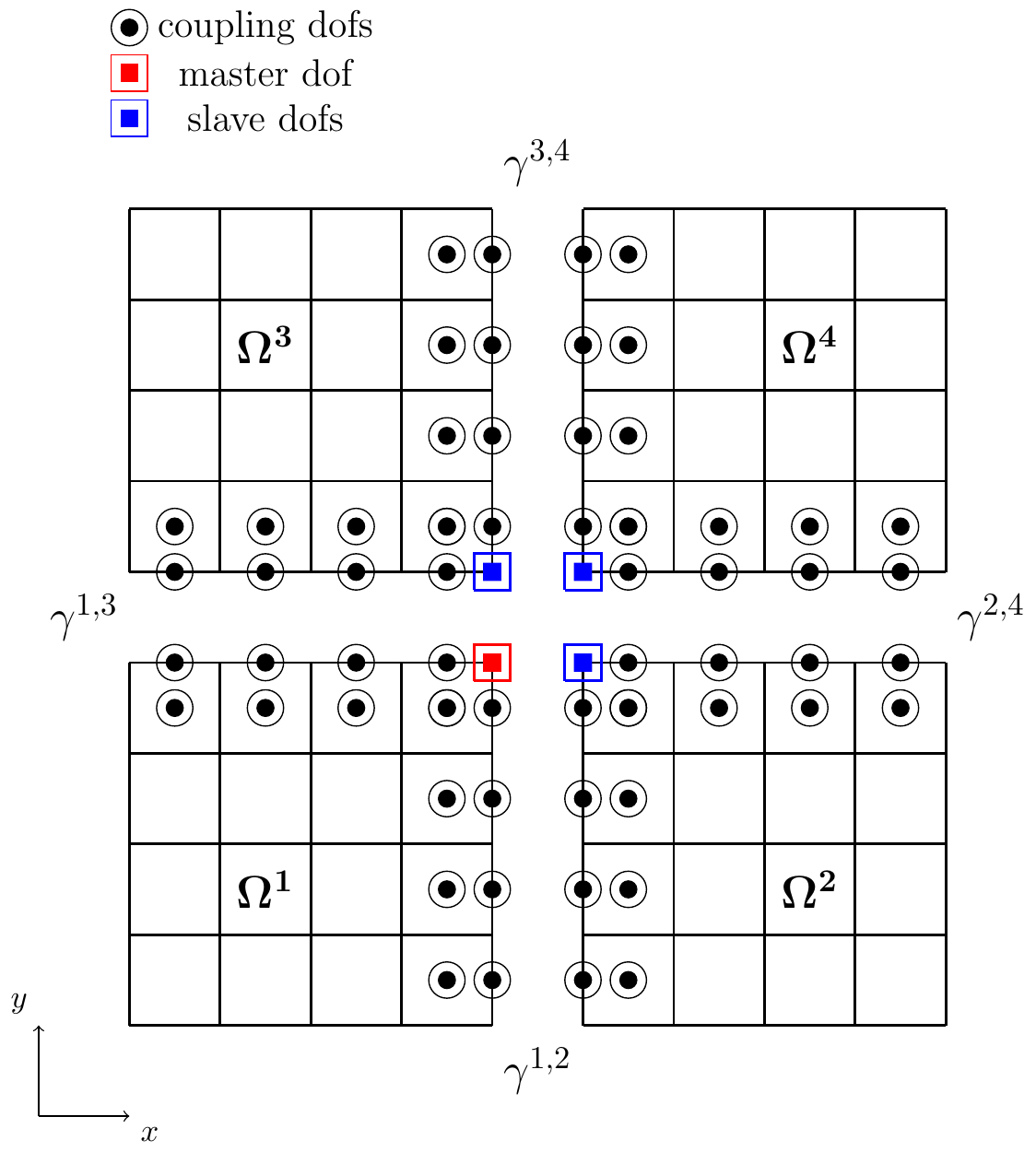}
	\caption{Cross-point modification, where the black dots represent the untouched control points associated to basis functions that give non-zero contribution to the interface coupling. The squares are the control variables used to impose the $C^0$ constraint, where we choose a master node (the red one) and the rest are labeled as slave nodes (the blue ones) and are eliminated from the system, see~\Cref{alg:apply_constrain}.}\label{fig:cross_point}
	\end{subfigure}
	\caption{Example of the dofs involved in the computation of the coupling integrals and cross-point modification in a four patches setup. For a correct interpretation of the colors, the reader is referred to the web version of this manuscript.}\label{fig:example_cp}
\end{figure}

In the literature of mortar methods, it is well-known that the treatment of cross-points requires extra considerations, see~\citep{Dittmann2020} and references therein for a discussion in the context of mortar coupling of isogeometric multi-patches. 
Analogously, our method also inherits the need for a cross-points modification. 
Indeed, in order to retain optimality of the method, a linear constraint must be imposed to the control variables meeting at the cross-point to ensure $C^0$-continuity. An example with four patches is depicted in~\Cref{fig:example_cp}, where in~\Cref{fig:dofs_interfaces} we depict the dofs associated to each coupling interface and in~\Cref{fig:cross_point} we visualize the imposition of the constraint.  
To explain the procedure, let us start from the following unconstrained system of equations:
\begin{align}\label{eq:reduced_problem_matrix}
\mathcal{A} \bm{u}_h = \bm{f} \, .
\end{align}
Now, the constraint can be incorporated easily into the standard linear system in a fully algebraic fashion, where a possible implementation is presented in~\Cref{alg:apply_constrain}. 

\begin{algorithm}[H] 
\begin{algorithmic}[1]
	\Procedure{Apply$\_C^0\_$constraint}{vector of dofs at cross-points $\bm{u}_{\text{cp}}$} 
		\State Label one dof in $\bm{u}_{\text{cp}}$ as master
		\State Label the remaining dofs in $\bm{u}_{\text{cp}}$ as slaves
		\State Build the rectangular matrix $\mathcal{C}$ representing the linear master-slaves constraints (see~\eqref{eq:build_constrain_matrix})
		\State Solve the reduced system $\widehat{\mathcal{A}} \widehat{\bm{u}}_h = \widehat{\bm{f}}$, where $\widehat{\mathcal{A}} = \mathcal{C}^\top \mathcal{A} \mathcal{C}$ and $\widehat{\bm{f}} = \mathcal{C}^\top \bm{f}$
		\State Recover the solution $\bm{u}_h$ from $\bm{u}_h = \mathcal{C} \widehat{\bm{u}}_h$
	\EndProcedure
\end{algorithmic} 
\caption{Algorithm for applying a $C^0$ constraint at a cross-point.}\label{alg:apply_constrain}
\end{algorithm}

The construction of the rectangular matrix $\mathcal{C}$ is best explained with an example. Let us assume that the dofs at the cross-point are numbered as $\bm{u}_{\text{cp}} = [{u}_{\text{cp} 1} \, {u}_{\text{cp} 2} \, {u}_{\text{cp} 3} \, {u}_{\text{cp} 4}]$. Now, without loss of generality, we pick ${u}_{\text{cp} 1}$ as the master control point and the rest as slave nodes. Then, the constraint can be expressed via the matrix $\mathcal{C}$ as follows:
\begin{align}\label{eq:build_constrain_matrix}
\bm{u}_h =
\begin{bmatrix}
u_1 \\
\vdots \\
{u}_{\text{cp} 1} \\
\vdots \\
{u}_{\text{cp} 2} \\
\vdots \\
{u}_{\text{cp} 3} \\
\vdots \\
{u}_{\text{cp} 4} \\
\vdots \\
u_{\text{ndof}}
\end{bmatrix} = 
\begin{blockarray}{cccccc}
& u_1 & \ldots & {u}_{\text{cp} 1} & \ldots & u_{\text{ndof}} \\
\begin{block}{c(ccccc)}
  u_1 & 1 & 0 & 0 & 0 & 0 \\
  \vdots & & \ddots & &  &  \\
  {u}_{\text{cp} 1} & 0 & 0 \ldots & 1 & 0 \ldots & 0 \\
  \vdots &  & \ddots  & &  &  \\
  {u}_{\text{cp} 2} & 0 & 0 \ldots & 1 & 0 \ldots & 0 \\
  \vdots &  &  & \ddots &  &  \\
  {u}_{\text{cp} 3} & 0 & 0 \ldots & 1 & 0 \ldots & 0 \\
  \vdots &  &  & \ddots &  &  \\
  {u}_{\text{cp} 4} & 0 & 0 \ldots & 1 & 0 \ldots & 0 \\
  \vdots &  &  & & \ddots &  \\
  u_{\text{ndof}} & 0 & 0 & 0 & 0 & 1 \\
\end{block}
\end{blockarray} \cdot  
\begin{bmatrix}
u_1 \\
\vdots \\
{u}_{\text{cp} 1} \\
\vdots \\
u_{\text{ndof}}
\end{bmatrix}
= \mathcal{C} \widehat{\bm{u}}_h \, ,
\end{align}
where ndof denotes the total number of degrees-of-freedom in the system. This procedure allows to eliminate the unknowns associated to the slave nodes from the system.

\section{A nested preconditioner based on the Schur Complement Reduction} \label{sec:preconditioner}

In this section, following the notation introduced in~\citep{Quarteroni20005} and building upon the work presented in~\citep{Liu2019,Liu2020} in the context of elastodynamics and hemodynamics, we present an efficient way to mitigate the detrimental effects on the condition number stemming from our choice of super-penalty parameters. This preconditioner is based on the approximate solution of the block factorization of the system matrix known as Schur Complement Reduction (SCR).
We remind the reader that before performing the algorithm described in the following, we apply a symmetric diagonal scaling to the system matrix. 

\subsection{The Schur Complement Reduction}

\begin{figure}
	\centering
	\includegraphics[width=0.6\textwidth]{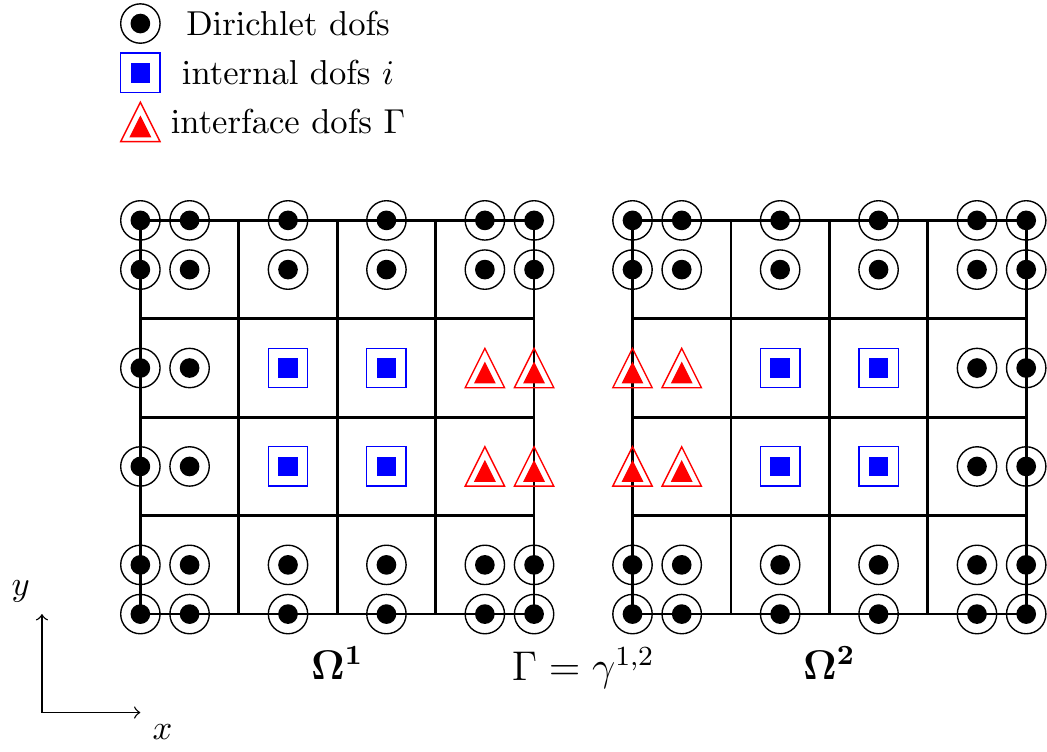}
	\caption{Example of reordering of the dofs in a two patches setup, discretized by B-splines of degree $p=2$, associated to the block system matrix $\mathcal{A}$.}\label{fig:system_dofs_reorder}
\end{figure}

We begin by reordering the matrix $\mathcal{A} \in \mathbb{R}^{\text{ndof} \times \text{ndof}}$ stemming from~\eqref{eq:reduced_problem} in blocks as follows:
\begin{align*}
\mathcal{A} =  
\begin{bmatrix}
\mathbf{A}_{i,i} & \mathbf{B}_{i,\Gamma} \\
\mathbf{B}_{i,\Gamma}^\top & \mathbf{C}_{\Gamma,\Gamma} \\
\end{bmatrix} \, ,
\end{align*}
where the subscripts $i$ and $\Gamma$ refer to internal and interface dofs, respectively, where an example is depicted in~\Cref{fig:system_dofs_reorder}. Let us remark that $\mathbf{A}_{i,i}$ is a block-diagonal matrix where every block is the matrix associated to an homogeneous Dirichlet problem (fully clamped) on the corresponding patch $\Omega^i$.
Moreover, with a slight abuse of notation, we assume that, if needed, $\mathcal{A}$ has already been modified to account for the constraints related to the cross-points introduced in the previous section.

Now, we can perform the following block factorization of $\mathcal{A}$:
\begin{align*}
\mathcal{A} = \mathcal{L} \mathcal{D} \mathcal{U} = 
\begin{bmatrix}
\mathbb{I} & \mathbf{0} \\
\mathbf{B}_{i,\Gamma}^\top \mathbf{A}_{i,i}^{-1} & \mathbb{I} \\
\end{bmatrix} 
\begin{bmatrix}
\mathbf{A}_{i,i} & \mathbf{0} \\
\mathbf{0} & \mathbf{S}_{\Gamma,\Gamma} \\
\end{bmatrix} 
\begin{bmatrix}
\mathbb{I} & \mathbf{A}_{i,i}^{-1} \mathbf{B}_{i,\Gamma} \\
\mathbf{0} & \mathbb{I} \\
\end{bmatrix} \, ,
\end{align*}
where we have introduced the Schur complement $\mathbf{S}_{\Gamma,\Gamma} := \mathbf{C}_{\Gamma,\Gamma} - \mathbf{B}_{i,\Gamma}^\top \mathbf{A}_{i,i}^{-1} \mathbf{B}_{i,\Gamma}$. 
Multiplying with $\mathcal{L}$ on both sides we get:
\begin{align}\label{eq:scr}
\begin{bmatrix}
\mathbf{A}_{i,i} & \mathbf{B}_{i,\Gamma} \\
\mathbf{0} & \mathbf{S}_{\Gamma,\Gamma} \\
\end{bmatrix}  
\begin{bmatrix}
\mathbf{x}_i \\
\mathbf{x}_\Gamma \\
\end{bmatrix} = 
\begin{bmatrix}
\mathbb{I} & \mathbf{0} \\
\mathbf{B}_{i,\Gamma}^\top \mathbf{A}_{i,i}^{-1} & \mathbb{I} \\
\end{bmatrix}^{-1}
\begin{bmatrix}
\mathbf{r}_i \\
\mathbf{r}_\Gamma \\
\end{bmatrix} = 
\begin{bmatrix}
\mathbb{I} & \mathbf{0} \\
- \mathbf{B}_{i,\Gamma}^\top \mathbf{A}_{i,i}^{-1} & \mathbb{I} \\
\end{bmatrix}
\begin{bmatrix}
\mathbf{r}_i \\
\mathbf{r}_\Gamma \\
\end{bmatrix} =
\begin{bmatrix}
\mathbf{r}_i \\
\mathbf{r}_\Gamma - \mathbf{B}_{i,\Gamma}^\top \mathbf{A}_{i,i}^{-1} \mathbf{r}_i \\
\end{bmatrix} \, .
\end{align}
We highlight that, up to this point, this factorization is performed in exact algebra.
Then, from~\eqref{eq:scr}, we can solve for $\mathbf{x}$ in a segregated fashion by exploiting~\Cref{alg:scr}. 
\begin{algorithm}[H] 
\begin{algorithmic}[1]
	\Procedure{Solution of $\mathcal{A} \mathbf{x} = \mathbf{r}$ based on SCR}{} 
	\State Solve for an intermediate solution $\mathbf{\hat{x}}_i$
	\begin{align}\label{eq:scr_alg_1}
		\mathbf{A}_{i,i} \mathbf{\hat{x}}_i = \mathbf{r}_i
	\end{align}	
	\State Update the interface residual $\mathbf{r}_\Gamma = \mathbf{r}_\Gamma - \mathbf{B}_{i,\Gamma}^\top \mathbf{\hat{x}}_i$	
	\State Solve for the interface solution $\mathbf{x}_\Gamma$ from the Schur equation
	\begin{align}\label{eq:scr_alg_2}
		\mathbf{S}_{\Gamma,\Gamma} \mathbf{x}_\Gamma = \mathbf{r}_\Gamma
	\end{align}	
	\State Update the internal residual $\mathbf{r}_i = \mathbf{r}_i - \mathbf{B}_{i,\Gamma} \mathbf{x}_\Gamma$	
	\State Solve for the internal solution $\mathbf{x}_i$ from
	\begin{align}\label{eq:scr_alg_3}
		\mathbf{A}_{i,i} \mathbf{x}_i = \mathbf{r}_i
	\end{align}		
	\EndProcedure
\end{algorithmic} 
\caption{SCR algorithm}\label{alg:scr}
\end{algorithm}    

Clearly, the Schur complement $\mathbf{S}_{\Gamma,\Gamma}$ is in practice expensive and often infeasible to compute explicitly. A way around this issue is given in Algorithm~\ref{alg:schur}, where we summarize a matrix-free procedure to apply the Schur complement to a vector.
\begin{algorithm}[H] 
\begin{algorithmic}[1]
	\Procedure{Application of $\mathbf{S}_{\Gamma,\Gamma}$ to a vector $\mathbf{x}_\Gamma$}{} 
	\State Compute the matrix-vector multiplication $\mathbf{\hat{x}}_\Gamma = \mathbf{C}_{\Gamma,\Gamma} \mathbf{x}_\Gamma$
	\State Compute the matrix-vector multiplication $\mathbf{\overline{x}}_\Gamma = \mathbf{B}_{i,\Gamma} \mathbf{x}_\Gamma$
	\State Solve for an intermediate solution $\mathbf{\tilde{x}}_\Gamma$ from
	\begin{align}\label{eq:schur_alg}
		\mathbf{A}_{i,i} \mathbf{\tilde{x}}_\Gamma = \mathbf{\overline{x}}_\Gamma
	\end{align}	
	\State Compute the matrix-vector multiplication $\mathbf{\overline{x}}_\Gamma = \mathbf{B}_{i,\Gamma}^\top \mathbf{\tilde{x}}_\Gamma$
	\State Return $\mathbf{\hat{x}}_\Gamma - \mathbf{\overline{x}}_\Gamma$
	\EndProcedure
\end{algorithmic} 
\caption{Algorithm for applying the Schur complement to a vector }\label{alg:schur}
\end{algorithm}    
\begin{remark}
As noted in~\citep{Liu2019}, the cost of the preconditioner is often dominated by the solution of the Schur system~\eqref{eq:scr_alg_2}. To reduce the computational burden of this step, we use as preconditioner a coarse approximation of the Schur complement obtained by applying only a few iterations of GMRES to $\mathbf{A}_{i,i}$ for assembling $\tilde{\mathbf{S}}_{\Gamma,\Gamma} = \mathbf{C}_{\Gamma,\Gamma} - \mathbf{B}_{i,\Gamma}^\top \tilde{\mathbf{A}}_{i,i}^{-1} \mathbf{B}_{i,\Gamma}$, where for efficiency we leverage again the FD algorithm. Although this choice works reasonably well for our numerical examples, we remark that more research is needed to find a robust (both in $h$ and $p$) and scalable preconditioner for the Schur complement and, more in general, for fourth-order PDEs.
\end{remark}

\subsection{Nested block preconditioner strategy based on SCR}
The main idea presented in~\citep{Liu2019} is to combine the robustness of the SCR factorization with the ease of application of a block preconditioners (such as SIMPLE or variants thereof~\citep{Quarteroni20005}). Indeed, we can build a preconditioner $\mathcal{P}_{\text{SCR}}$ based on an approximate factorization of~\eqref{eq:scr}, where~\Cref{eq:scr_alg_1,eq:scr_alg_2,eq:scr_alg_3} are solved within a prescribed tolerance. 
Given that $\mathcal{P}_{\text{SCR}}$ changes its algebraic definition at every iteration, following~\citep{Liu2019}, we employ a flexible GMRES algorithm (FGMRES) as the iterative method for the most outer solve $\mathcal{A} \mathbf{x} = \mathbf{r}$.
At each iteration of FGMRES, we can apply the preconditioner $\mathcal{P}_{\text{SCR}}$ via~\Cref{alg:scr}, where this entails the solution of the blocks $\mathbf{A}_{i,i}$ and $\mathbf{S}_{\Gamma,\Gamma}$. This part of the algorithm is denoted as intermediate solver. Last, since we do not assemble the Schur complement explicitly, but we apply its action on a vector through~\Cref{alg:schur}, we perform a final solve for $\mathbf{A}_{i,i}$ in~\eqref{eq:schur_alg}, denoted as inner solver.
The final performance of the preconditioner is therefore determined by the prescribed tolerances for the outer, intermediate and inner layers, respectively, where the objective is finding a good balance between the computational cost and the robustness of the method. In the following, we denote the aforementioned tolerances by $\eta_o, \eta_t$ and $\eta_n$ for the outer, intermediate and inner layers, respectively.

\subsubsection{A preconditioner based of the Fast Diagonalization (FD) algorithm}

Since each outer iteration of the nested preconditioner is based on the solution of three systems involving the block $\mathbf{A}_{i,i}$, an efficient and robust preconditioner for this block is required. In this work, we extend the isogeometric preconditioner studied in~\citep{Tani2016,Montardini2018}, based on the Fast Diagonalization algorithm, to the Kirchhoff plate problem.
In the following, we focus our derivation on the single-patch case. The extension to the multi-patch case is straightforward by construction, since the block $\mathbf{A}_{i,i}$ is formed by disjoint sub-blocks associated to each patch $\Omega^i$.


Now, exploiting the tensor product structure of the B-spline basis at the patch level, let us introduce the preconditioner $\mathcal{P}_{\text{FD}}$ in Kronecker form as:
\begin{align}\label{eq:precond_FD}
\mathcal{P}_{\text{FD}} = M_1 \otimes K_2 + K_1 \otimes M_2 \, ,
\end{align}
where $M_k$ and $K_k$ with $k=1,2$ refer to the one-dimensional, parametric mass and hessian matrices associated to the $k$-th parametric dimension, respectively. They can be expanded as follows:
\begin{align*}
\left[ M_k \right]_{i,j} &= \int_0^1 b_{i,p} (\eta_k) \, b_{j,p} (\eta_k) \text{ d}\eta_k  \\
\left[ K_k \right]_{i,j} &= \int_0^1 b_{i,p}^{\prime\prime} (\eta_k) \, b_{j,p}^{\prime\prime} (\eta_k) \text{ d}\eta_k \, ,
\end{align*}
where $b$ indicates the univariate B-spline basis functions introduced in~\Cref{sec:bsplines}.
Then, analogously to~\citep{Montardini2018spacetime}, we partially include the geometry and physical coefficients inside the preconditioner.
In particular, let us denote by $\mathfrak{C}$ the following function:
\begin{align*}
\mathfrak{C}(\bm{\eta}) = D \left( \big| \big| J_{\mathbf{F}}^{-1} \big| \big|_2 \right)^4 | \text{det }(J_{\mathbf{F}}) | \, ,
\end{align*} 
where we recall that $J_{\mathbf{F}}$ represents the jacobian of the B-spline parametrization $\mathbf{F}$ and $D$ is the flexural stiffness of the plate.
Now, as explained in~\citep[Appendix A.3]{Montardini2018spacetime}, we perform a separation of variables on $\mathfrak{C}$ such that we can write:
\begin{align*}
\mathfrak{C}(\bm{\eta}) \approx \widetilde{\mathfrak{C}}(\bm{\eta}) =
\begin{bmatrix}
\omega_1(\eta_1) \tau_2(\eta_2) & 0 \\
0 & \tau_1(\eta_1) \omega_2(\eta_2) 
\end{bmatrix} \, , 
\end{align*}
where this matrix is evaluated at each quadrature point.
With this, we can modify the preconditioner given in~\eqref{eq:precond_FD} to partially account for the geometry and coefficients information as follows:
\begin{align}\label{eq:precond_FD_geo}
\mathcal{P}_{\text{FD}}^\mathbf{F} = \widetilde{M}_1 \otimes \widetilde{K}_2 + \widetilde{K}_1 \otimes \widetilde{M}_2 \, ,
\end{align}
where
\begin{align*}
\left[ \widetilde{M}_k \right]_{i,j} &= \int_0^1 \omega_k(\eta_k) b_{i,p} (\eta_k) \, b_{j,p} (\eta_k) \text{ d}\eta_k  \\
\left[ \widetilde{K}_k \right]_{i,j} &= \int_0^1 \tau_k(\eta_k) b_{i,p}^{\prime\prime} (\eta_k) \, b_{j,p}^{\prime\prime} (\eta_k) \text{ d}\eta_k \, .
\end{align*}
Finally, each iteration of the iterative solver requires the solution of the following system:
\begin{align}\label{eq:residual_FD}
\mathcal{P}_{\text{FD}}^\mathbf{F} s = r \, ,
\end{align}
where $r$ denotes the current residual. Due to the tensor structure of the preconditioner, we can rewrite~\eqref{eq:residual_FD} as a Sylvester matrix equation~\citep{Simoncini2016}:
\begin{align*}
\widetilde{M}_2 S  \widetilde{K}_1 +  \widetilde{K}_2 S  \widetilde{M}_1 = R \, ,
\end{align*}
where $s = \text{vec}(S)$ and $r = \text{vec}(R)$.
\begin{remark}
Let us recall that for any matrix $Z \in \mathbb{R}^{r \times c}$ the operator $\emph{vec}(Z)$ gives as output the vector $z \in \mathbb{R}^{r c}$ formed by stacking the columns of $Z$.  
\end{remark}
\noindent Let us now consider the generalized eigendecomposition of the matrix pencils $(\widetilde{K}_1, \widetilde{M}_1)$ and $(\widetilde{K}_2, \widetilde{M}_2)$, respectively, as:
\begin{align}\label{eq:generalized_eigendecomposition}
\widetilde{K}_1 U_1 &= \widetilde{M}_1 U_1 D_1  \nonumber \\
\widetilde{K}_2 U_2 &= \widetilde{M}_2 U_2 D_2 \, .
\end{align}
Here, $D_1$ and $D_2$ are diagonal matrices containing the eigenvalues of $\widetilde{M}_1^{-1} \widetilde{K}_1$ and $\widetilde{M}_2^{-1} \widetilde{K}_2$, respectively. Further, $U_1$ and $U_2$ are defined as:
\begin{align*}
U_1^\top \widetilde{M}_1 U_1 &= \mathbb{I}  \\ 
U_2^\top \widetilde{M}_2 U_2 &= \mathbb{I} \, . 
\end{align*}
With these definitions at hand, we can rewrite~\eqref{eq:precond_FD_geo} in Kronecker form as:
\begin{align*}
\left( U_1 \otimes U_2 \right)^{-\top} \left( D_1 \otimes \mathbb{I} + \mathbb{I} \otimes D_2 \right) \left( U_1 \otimes U_2 \right)^{-1} s = r \, ,
\end{align*}
where the preconditioner can be efficiently applied via~\Cref{alg:FD}.
\begin{algorithm}[H] 
\begin{algorithmic}[1]
	\Procedure{Update of the iteration residual via the FD method}{} 
	\State Compute the generalized eigendecomposition in~\eqref{eq:generalized_eigendecomposition}
	\State Compute the intermediate result $\tilde{r} = \left( U_1 \otimes U_2 \right)^{\top} r$
	\State Compute the intermediate residual $\tilde{s} = \left( D_1 \otimes \mathbb{I} + \mathbb{I} \otimes D_2 \right)^{-1} \tilde{r}$
	\State Return $s = \left( U_1 \otimes U_2 \right) \tilde{s}$
	\EndProcedure
\end{algorithmic} 
\caption{FD method for applying $\mathcal{P}_{\text{FD}}^\mathbf{F}$  }\label{alg:FD}
\end{algorithm} 

\begin{remark}
We remark that the application of the nested preconditioner $\mathcal{P}_{\text{SCR}}$ combined with $\mathcal{P}_{\text{FD}}^\mathbf{F}$ can be implemented in a fully matrix-free framework. Furthermore, although not investigated in this work, the patch-wise block structure of $\mathbf{A}_{i,i}$ could be further exploited for parallelization. 
\end{remark}  

\noindent For the sake of conciseness, we do not provide here further details of the FD algorithm, but we refer to~\citep{Tani2016,Montardini2018} for a thorough theoretical and numerical investigation of the method in the scope of isogeometric analysis.

\renewcommand{\graphDir}{pictures/numericalExamples/graphs}
\renewcommand{\dataDir}{pictures/numericalExamples/data}
\renewcommand{\picsDir}{pictures/numericalExamples/pics}

\section{Numerical Examples}
\label{sec:numericalExamples}

In this section we assess the performance of the proposed coupling method with several numerical examples defined on multi-patch geometries. 
All the numerical experiments presented in the following have been implemented in the open-source and free Octave/Matlab package \textit{GeoPDEs}~\citep{Vazquez2016}, a software designed for the solution of partial differential equations in the context of isogeometric analysis.  	


\subsection{A four patches example with non-matching curved interfaces}

\begin{figure}
	\centering
	\begin{subfigure}[t]{0.495\textwidth}
		\centering
	\includegraphics[width=0.95\textwidth]{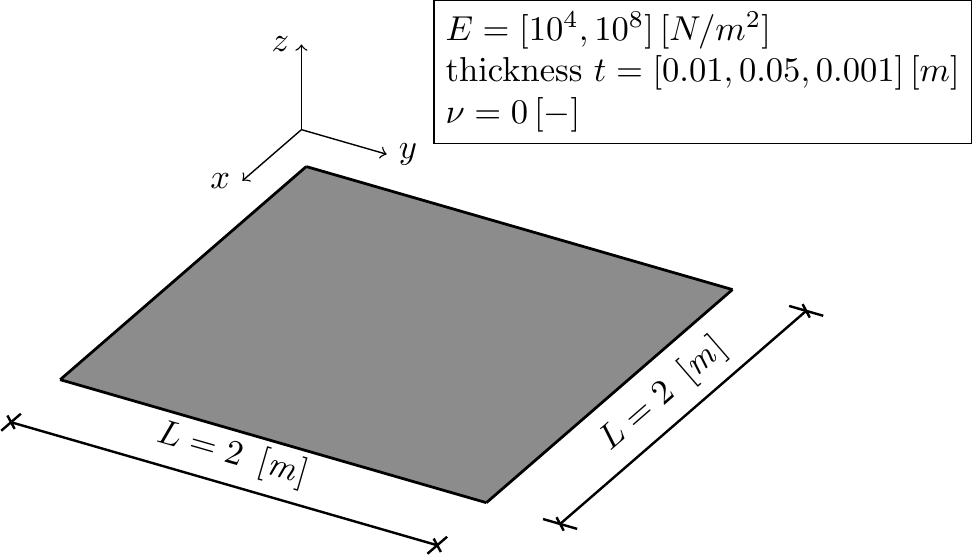}
		\caption{Geometry setup and physical parameters.}
	\end{subfigure}
	\hfill
	\begin{subfigure}[t]{0.495\textwidth}
		\centering
	\includegraphics[width=0.85\textwidth]{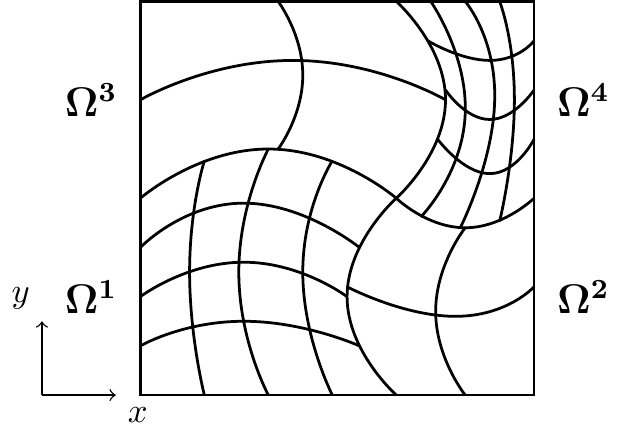}
		\caption{Initial discretization.}
	\end{subfigure}	
	\caption{Problem setup and initial multi-patch non-conforming discretization for the curved four patches example.}
\label{fig:setup_curved_4patches}
\end{figure}

\begin{figure}
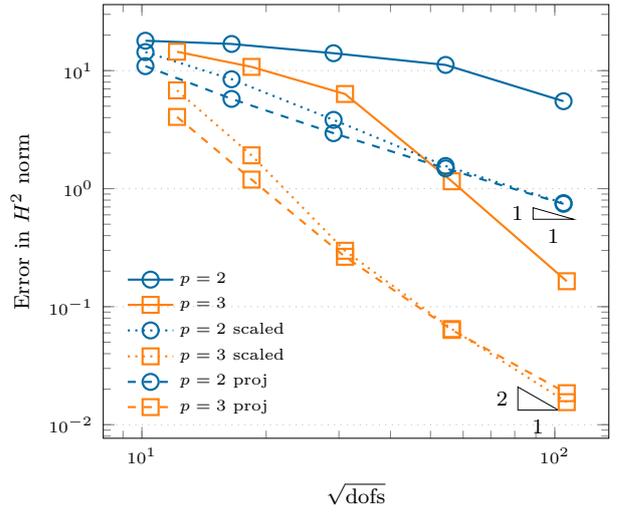
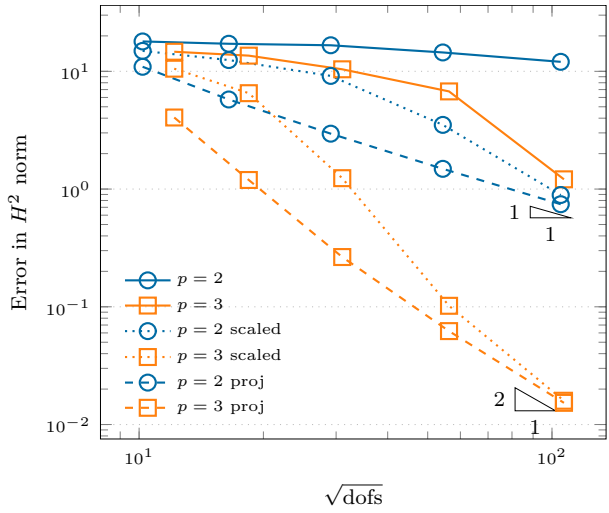
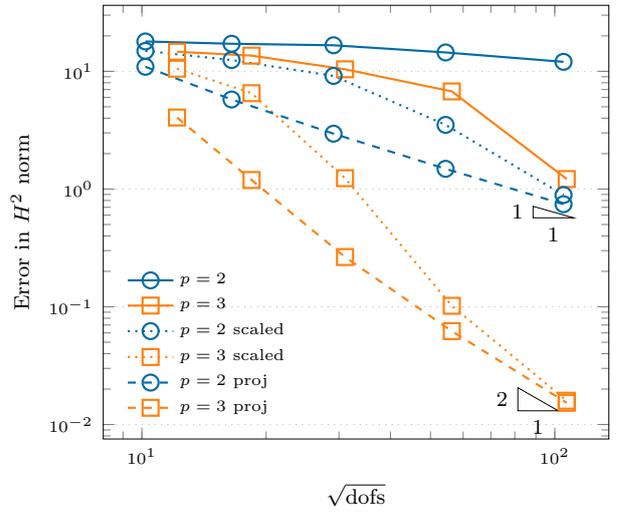
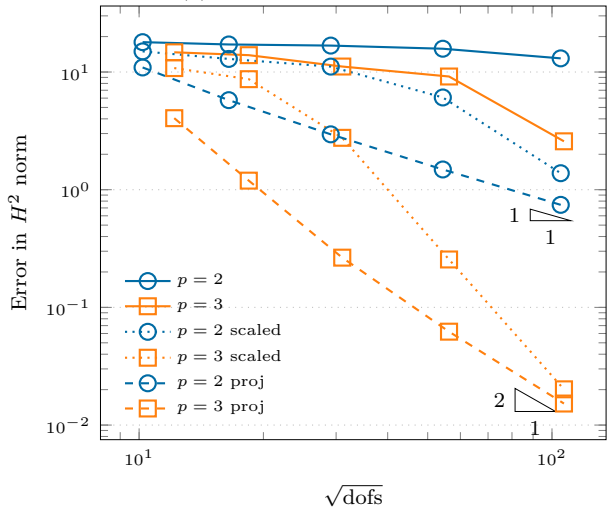
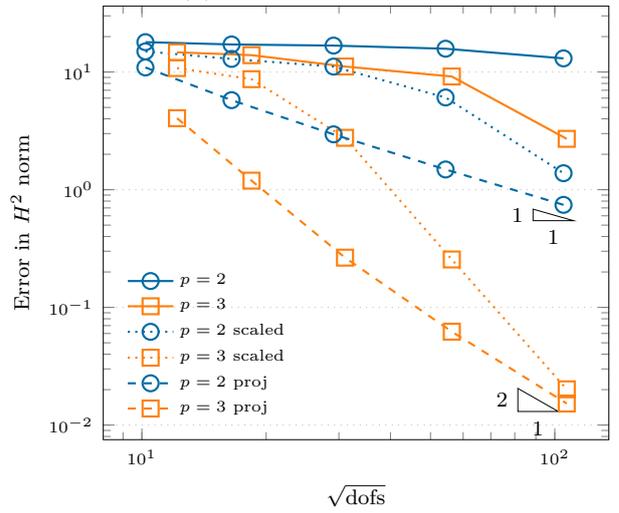

	\centering
	\begin{subfigure}[t]{0.495\textwidth}
		\centering
	\input{\graphDir/convergence_bilaplace_curved_4patches_E1e4_t005.tex}
		\caption{$E = 10^4, t = 0.05$.}
	\end{subfigure}
	\hfill
	\begin{subfigure}[t]{0.495\textwidth}
		\centering
	\input{\graphDir/convergence_bilaplace_curved_4patches_E1e8_t005.tex}
		\caption{$E = 10^8, t = 0.05$.}
	\end{subfigure}	
	\begin{subfigure}[t]{0.495\textwidth}
		\centering
	\input{\graphDir/convergence_bilaplace_curved_4patches_E1e4_t001.tex}
		\caption{$E = 10^4, t = 0.01$.}
	\end{subfigure}
	\hfill
	\begin{subfigure}[t]{0.495\textwidth}
		\centering
	\input{\graphDir/convergence_bilaplace_curved_4patches_E1e8_t001.tex}
		\caption{$E = 10^8, t = 0.01$.}
	\end{subfigure}		
	\begin{subfigure}[t]{0.495\textwidth}
		\centering
	\input{\graphDir/convergence_bilaplace_curved_4patches_E1e4_t0005.tex}
		\caption{$E = 10^4, t = 0.005$.}
	\end{subfigure}
	\hfill
	\begin{subfigure}[t]{0.495\textwidth}
		\centering
	\input{\graphDir/convergence_bilaplace_curved_4patches_E1e8_t0005.tex}
		\caption{$E = 10^8, t = 0.005$.}
	\end{subfigure}		
	\caption{Convergence study of the error measured in the $H^2$ norm in the non-matching case for four patches with curved interface example for different Young moduli and values of the thickness, B-splines of degree $p=2,3$. Comparison of a classic penalty method, the scaled version with respect to the problem parameters proposed in~\citep{Herrema2019} (\textit{scaled}) and our projection approach (\textit{proj}).}
\label{fig:convergence_bilaplace_curved_4patches_parameters}
\end{figure}

In this example we consider the computational domain $\Omega = [0,2] \times [0,2]$ depicted in~\Cref{fig:setup_curved_4patches}, split into four subdomains $\Omega^i$. We remark that all meshes are non-conforming at every coupling interface, as the irrational factor $\sqrt{2}/100$ has been used to shift the interface knots. 
The body source and boundary data are computed such that the exact solution is smooth and it reads:
\begin{align*}
u^{\text{ex}} = \sin(\pi x) \cos(\pi x)  \, .
\end{align*}
This setup is used to test the robustness of our method in the case of severe non-matching discretizations and with respect to the problem parameters.
To this end,  we present the convergence results for all combinations of Young's moduli $E = [10^4, 10^8] \, [Pa]$ and thickness of the plate $t = [0.05, 0.01, 0.005] \, [m]$, where we set the Poisson's ratio $\nu = 0 \, [-]$. We compare our method to a classical penalty approach, where we set $\alpha_{\text{defl}}^\ell = \alpha_{\text{rot}}^\ell = 10^4 E, \ell = 1,\ldots,L$, and to a choice of penalty parameters scaled with respect to the physical parameters as proposed in~\citep{Herrema2019}. In particular, they read:
\begin{align*}
\alpha_{\text{defl}}^\ell &= \delta \frac{E t}{h_\ell(1 - \nu^2)} \\	
\alpha_{\text{rot}}^\ell &= \delta \frac{E t^3}{12 h_\ell(1 - \nu^2)} \, ,
\end{align*}
where the user-defined parameter $\delta = 10^3$ is chosen. From the results in~\Cref{fig:convergence_bilaplace_curved_4patches_parameters}, we observe that the projection strategy shows robustness with respect to the input parameters and allows for an easy treatment of locking phenomena, where optimal convergence rates are attained also for very coarse meshes.

\begin{figure}
	\centering
	\begin{subfigure}[t]{0.495\textwidth}
		\centering
	\input{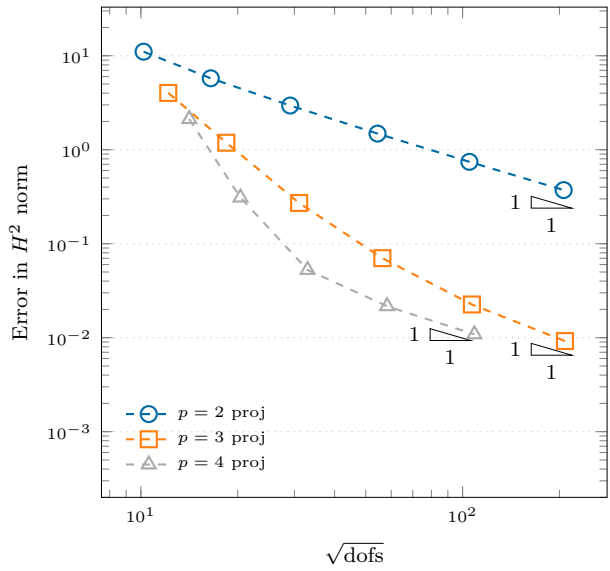}
		\caption{Without $C^0$ constraint.}
	\end{subfigure}
	\hfill
	\begin{subfigure}[t]{0.495\textwidth}
		\centering
	\input{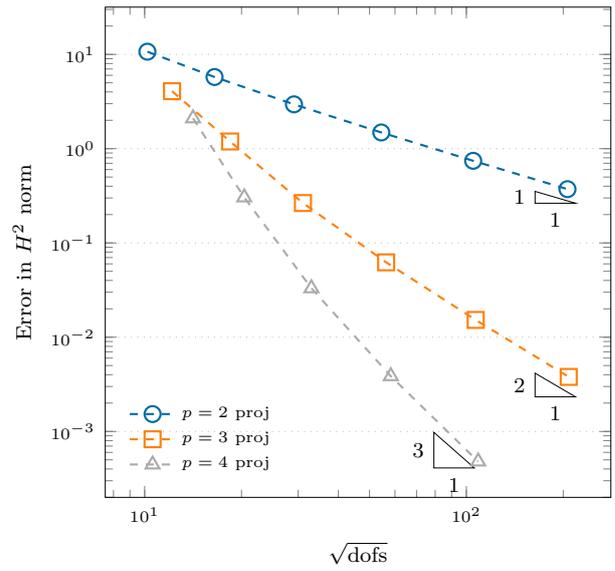}
		\caption{With $C^0$ constraint.}
	\end{subfigure}	
	\caption{Convergence study of the error in the $H^2$ norm in the non-matching case for the curved four patches example. Influence of imposing a $C^0$ constraint at the cross point.}
\label{fig:convergence_bilaplace_curved_4patches}
\end{figure}

\begin{figure}
	\centering
	\begin{subfigure}[t]{0.495\textwidth}
		\centering
	\includegraphics[width=0.995\textwidth]{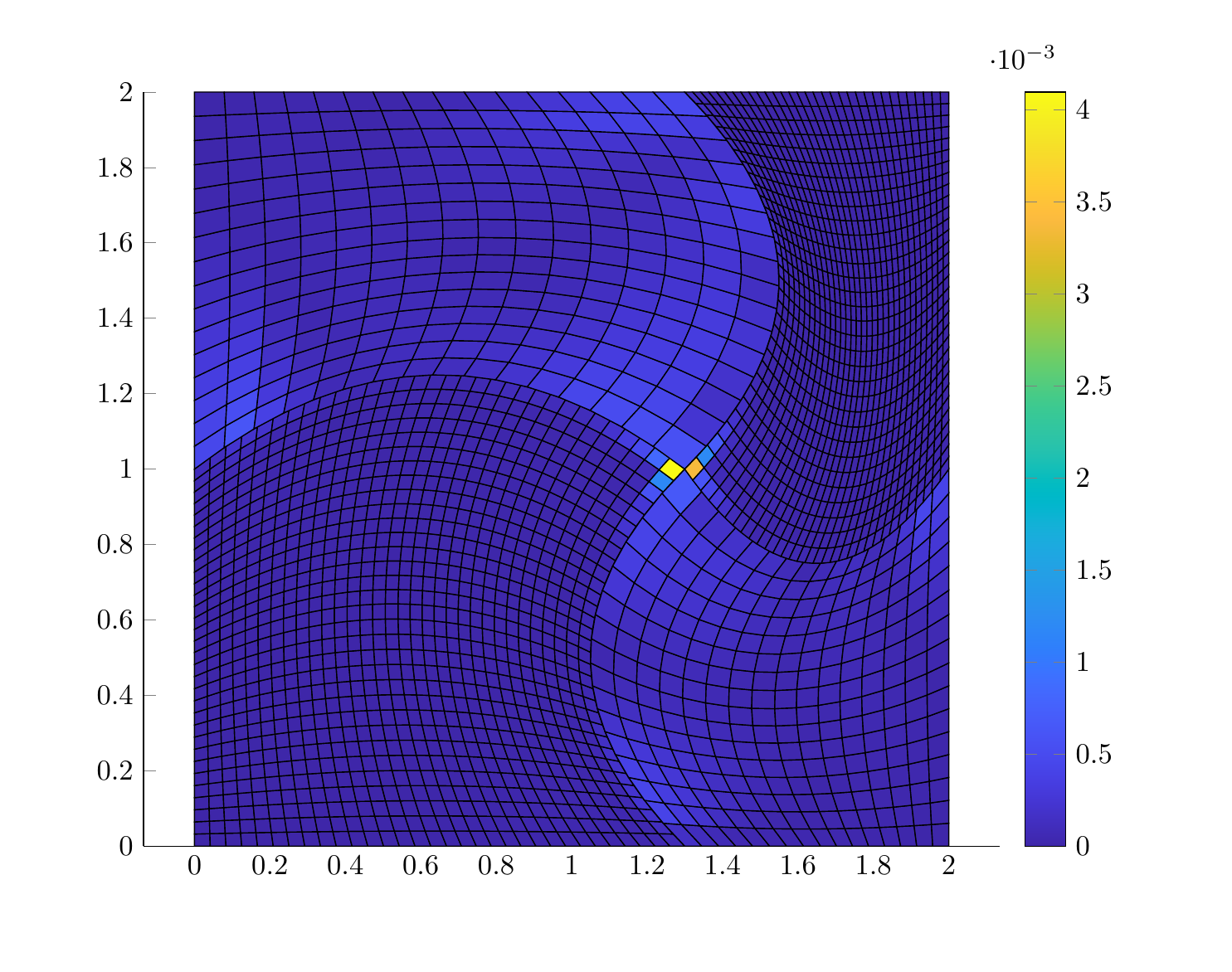}
		\caption{Without $C^0$ constraint.}
	\end{subfigure}
	\hfill
	\begin{subfigure}[t]{0.495\textwidth}
		\centering
	\includegraphics[width=0.995\textwidth]{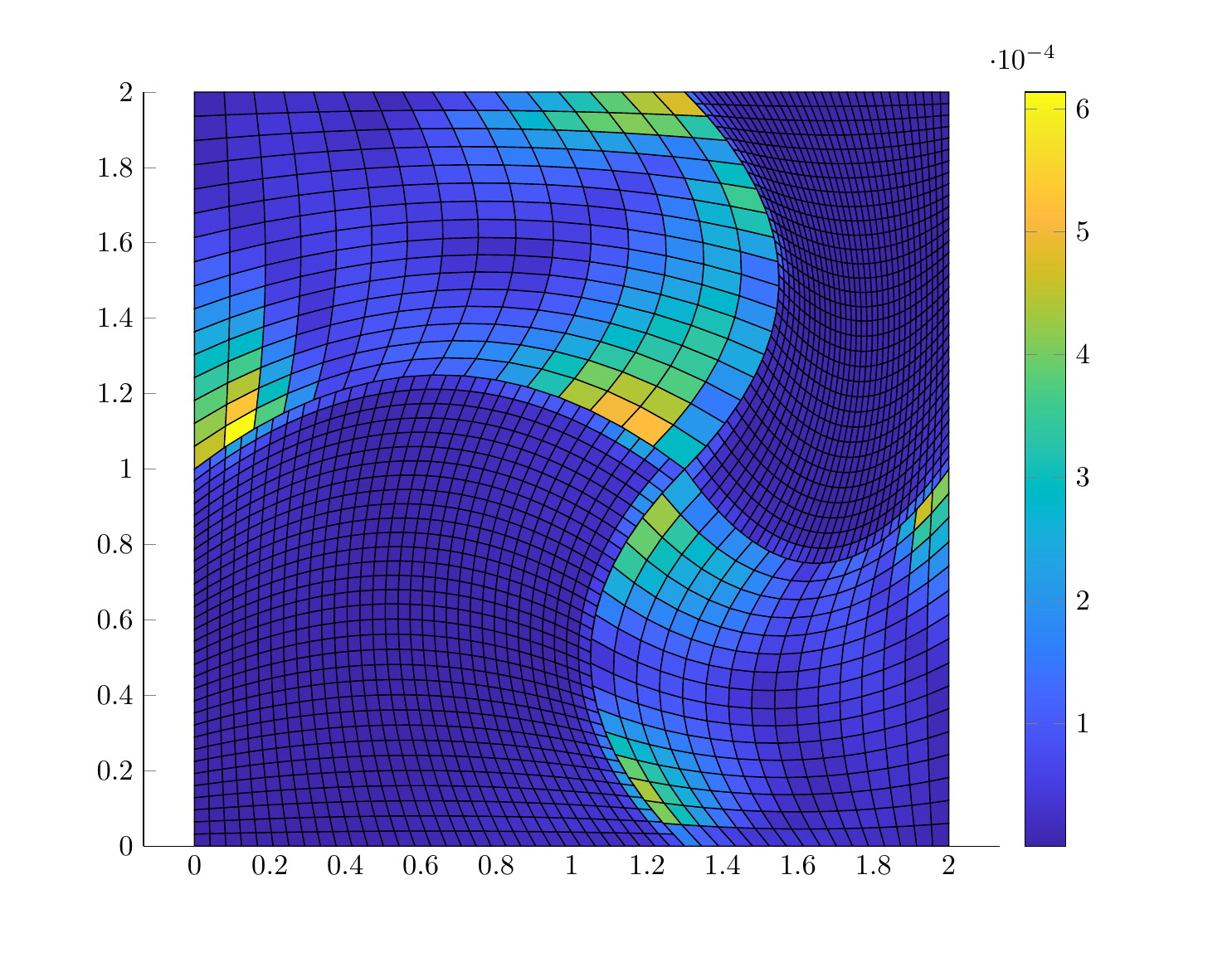}
		\caption{With $C^0$ constraint.}
	\end{subfigure}	
	\caption{Element-wise plot of the error in the $H^2$ norm in the non-matching case for the curved four patches example, B-splines of degree $p=4$. Influence of imposing a $C^0$ constraint at the cross point, notice the difference of one order of magnitude used in the two colorbars.}
\label{fig:error_plot_bilaplace_curved_4patches}
\end{figure}

In~\Cref{fig:convergence_bilaplace_curved_4patches} the convergence behavior of the error measured in the $H^2$ norm with and without the imposition of the $C^0$ constraint at the cross-point is plotted. We observe that the loss of accuracy hinders the convergence for $p=3,4$, whereas the expected optimal rates of convergence are recovered in all cases when the linear constraint is imposed to the system. This is further highlighted in~\Cref{fig:error_plot_bilaplace_curved_4patches}, where the element-wise $H^2$ error is depicted for a discretization of degree $p=4$, without and with the constraint, respectively. On one hand, we remark how the error is concentrated and much higher in the elements around the cross-point, spoiling the optimal convergence, when the constraint is not imposed. On the other hand, with the linear constraint, we recover optimal convergence properties of the method. 

Finally, for this example we also analyze the performance of the nested preconditioner. In~\Cref{tab:precond} we report the iterations needed by the external solver and in brackets the average number of intermediate iterations, for several degrees of the discretization $p=2,3$, and we compare it with a classical diagonally preconditioned conjugate gradient (PCG), a PCG where an incomplete LU (ILU) is used as preconditioner and a GMRES preconditioned with ILU. All the results refer to a global tolerance $\eta_o$ of $10^{-10}$ and, for the nested SCR-FGMRES strategy, the intermediate and inner tolerances $\eta_t$ and $\eta_n$ are set to $10^{-6}$. Further, the Schur complement is preconditioned by an approximation $\tilde{\mathbf{S}}_{\Gamma,\Gamma}$ obtained with a maximum of 6 iterations of GMRES.
For the sake of completeness, we perform the same test with the choice of penalty parameters studied in~\citep{Herrema2019}. The results are summarized in~\Cref{tab:precond_josef}, where we observe no substantial difference regarding the iterations needed to solve the system compared to the case where our choice of parameters is employed. This suggests that the proposed preconditioner is robust with respect to the penalty factors and it is also suitable to precondition systems stemming from other penalty approaches.

\begin{table}[h]
\begin{subtable}{\textwidth}
\begin{center}
    \begin{tabular}{| l | c | c | c | c |}
    \hline
     & 256 el. & 1024 el. & 4096 el. & 16384 el.  \\ \hline
    Diagonally scaled PCG & 792 & 953 & -- & -- \\ \hline
    PCG with ILU &  111 & 980 & -- & -- \\ \hline
    GMRES with ILU & 63 & 174 & 402 & -- \\ \hline
    Nested SCR-FGMRES &  3 (21.6/3.3/20.3) & 3 (36/6/30.3) & 4 (51/17.5/40.7) & 6 (66.5/51/51.3) \\ \hline
    \hline
    \end{tabular}
    \caption{$p=2$.} \label{tab:precond_p2}
\end{center} 
\end{subtable}
\begin{subtable}{\textwidth}
\begin{center}
    \begin{tabular}{| l | c | c | c | c |}
    \hline
     & 256 el. & 1024 el. & 4096 el. & 16384 el.  \\ \hline
    Diagonally scaled PCG & 921 & -- & -- & -- \\ \hline
    PCG with ILU &  53 & 221 & -- & -- \\ \hline
    GMRES with ILU &  35 & 73 & -- & -- \\ \hline
    Nested SCR-FGMRES & 3 (26.6/4/24) & 3 (41.6/9/35) & 4 (58/25.5/45.7) & 6 (76.3/68.3/55.3) \\ \hline
    \hline
    \end{tabular}
   \caption{$p=3$.} \label{tab:precond_p3}
\end{center} 
\end{subtable}
   \caption{Number of iterations needed by different iterative methods, $p=2,3$, as a function of the elements (el.). For the nested SCR-FGMRES, the numbers in brackets indicate the average number of intermediate iterations needed to solve~\Cref{eq:scr_alg_1,eq:scr_alg_2,eq:scr_alg_3} in~\Cref{alg:scr}, respectively. Iterations marked with -- did not reached convergence within the prescribed 1000 maximum number of iterations.} \label{tab:precond}
\end{table}  

\begin{table}[h]
\begin{subtable}{\textwidth}
\begin{center}
    \begin{tabular}{| l | c | c | c | c |}
    \hline
     & 256 el. & 1024 el. & 4096 el. & 16384 el.  \\ \hline
    Diagonally scaled PCG & 992 & -- & -- & -- \\ \hline
    PCG with ILU &  175 & -- & -- & -- \\ \hline
    GMRES with ILU & 79 & 209 & 478 & -- \\ \hline
    Nested SCR-FGMRES &  3 (22.3/4/21.6) & 3 (32/7.3/28.3) & 4 (44.5/20.2/34.7) & 5 (63.2/55/50.2) \\ \hline
    \hline
    \end{tabular}
    \caption{$p=2$.} \label{tab:precond_p2_josef}
\end{center} 
\end{subtable}
\begin{subtable}{\textwidth}
\begin{center}
    \begin{tabular}{| l | c | c | c | c |}
    \hline
     & 256 el. & 1024 el. & 4096 el. & 16384 el.  \\ \hline
    Diagonally scaled PCG & 778 & -- & -- & -- \\ \hline
    PCG with ILU &  697 & -- & -- & -- \\ \hline
    GMRES with ILU & 90 & 209 & -- & -- \\ \hline
    Nested SCR-FGMRES &  3 (28.6/5.3/27.3) & 3 (39.3/13/34.6) & 5 (55.8/29.4/41) & 6 (74/73.3/48.3) \\ \hline
    \hline
    \end{tabular}
   \caption{$p=3$.} \label{tab:precond_p3_josef}
\end{center} 
\end{subtable}
\caption{Number of iterations needed by different iterative methods, $p=2,3$, as a function of the elements (el.) for the parameters proposed in~\citep{Herrema2019}. For the nested SCR-FGMRES, the numbers in brackets indicate the average number of intermediate iterations needed to solve~\Cref{eq:scr_alg_1,eq:scr_alg_2,eq:scr_alg_3} in~\Cref{alg:scr}, respectively. Iterations marked with -- did not reached convergence within the prescribed 1000 maximum number of iterations.} \label{tab:precond_josef}
\end{table}  

In~\Cref{tab:precond_tol} we study the influence of the intermediate and inner tolerances on the number of outer iterations required by the FGMRES solver, on a fixed mesh of 4096 elements, for B-splines of degree $p=2,3$. We note that as the chosen tolerances become smaller and smaller, we recover the algebraically exact SCR method, where in the limit the algorithm converges in one iteration. We also remark that finding an optimal choice for these parameters is, to the best of the authors' knowledge, still an open question in the community. 

\begin{table}[h]
\begin{center}
    \begin{tabular}{| l | c | c | c | c | c |}
    \hline
     & $\eta_t=\eta_n=10^{-4}$ & $\eta_t=\eta_n=10^{-5}$ & $\eta_t=\eta_n=10^{-6}$ & $\eta_t=\eta_n=10^{-8}$ & $\eta_t=\eta_n=10^{-10}$ \\ \hline
    $p=2$ & 11 & 5 & 4 & 3 & 2 \\ \hline
    $p=3$ & 13 & 7 & 4 & 3 & 2 \\ \hline
    \hline
    \end{tabular}
    \caption{Influence of the intermediate and inner tolerances $\eta_t$ and $\eta_n$ (where we always set $\eta_t = \eta_n$) on the number of outer iterations needed by the FGMRES solver, $p=2,3$, on a fixed mesh with 4096 elements.} \label{tab:precond_tol}
\end{center} 
\end{table}  

\subsection{A nine patches geometry}

\begin{figure}
	\centering
	\begin{subfigure}[t]{0.495\textwidth}
		\centering
	\includegraphics[width=0.95\textwidth]{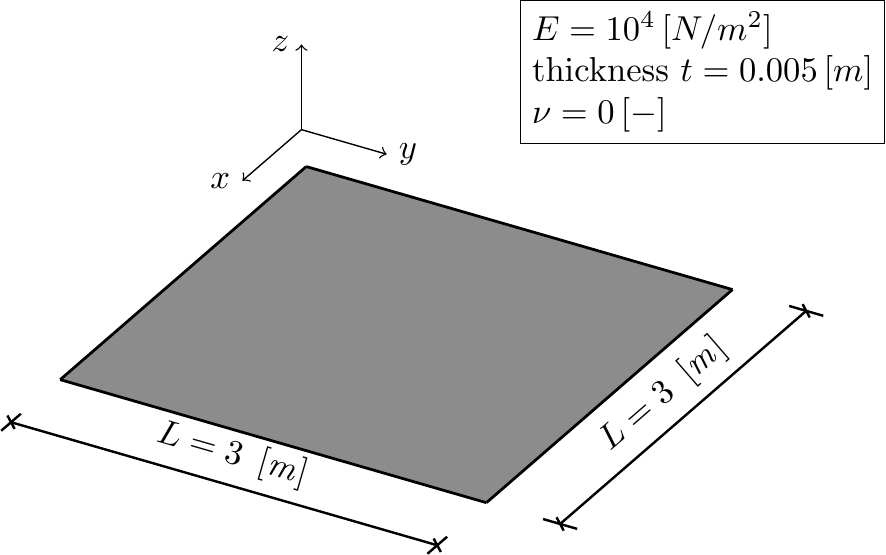}
		\caption{Geometry setup and physical parameters.}
	\end{subfigure}
	\hfill
	\begin{subfigure}[t]{0.495\textwidth}
		\centering
	\includegraphics[width=0.75\textwidth]{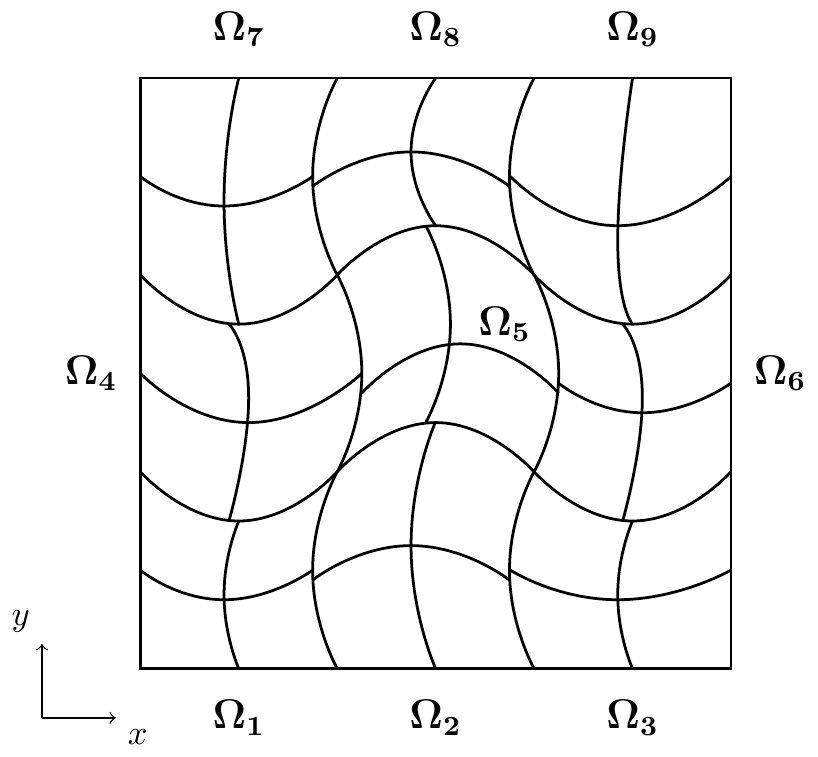}
		\caption{Initial discretization.}
	\end{subfigure}	
	\caption{Problem setup and initial multi-patch non-conforming discretization for the nine patches example.}
\label{fig:setup_9patches}
\end{figure}

\begin{figure}
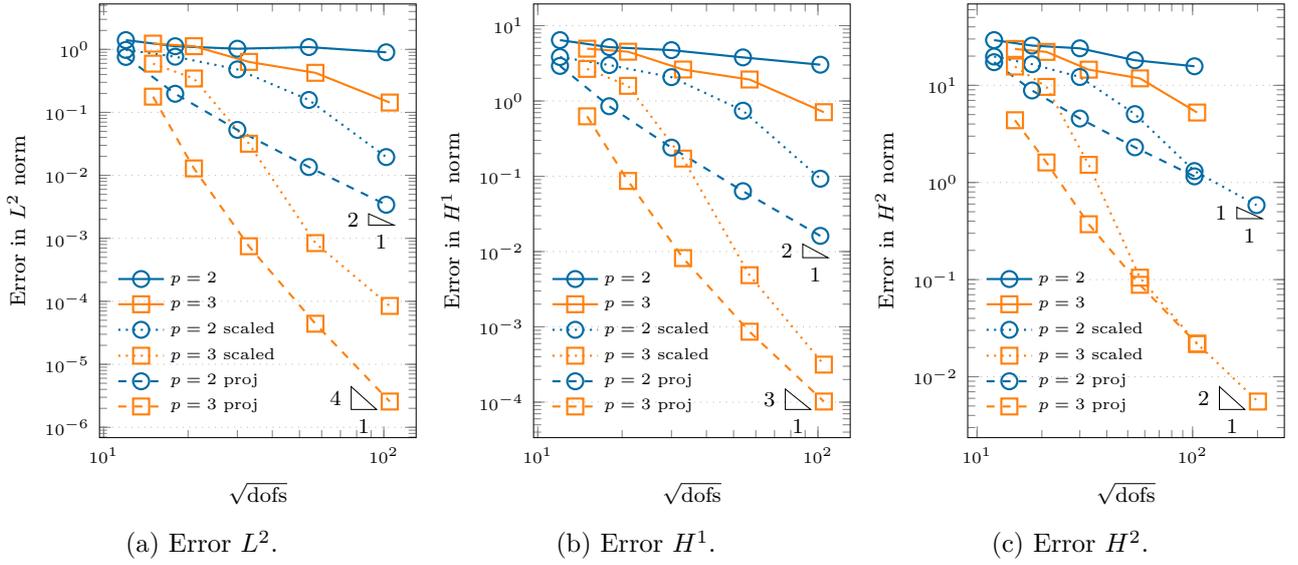

	\centering
	\begin{subfigure}[t]{0.32\textwidth}
		\centering
	\input{\graphDir/convergence_bilaplace_9patches_nonconforming_l2.tex}
		\caption{Error $L^2$.}
	\end{subfigure}
	\hfill
	\begin{subfigure}[t]{0.32\textwidth}
		\centering
	\input{\graphDir/convergence_bilaplace_9patches_nonconforming_h1.tex}
		\caption{Error $H^1$.}
	\end{subfigure}
	\hfill
	\begin{subfigure}[t]{0.32\textwidth}
		\centering
	\input{\graphDir/convergence_bilaplace_9patches_nonconforming_h2.tex}
		\caption{Error $H^2$.}
	\end{subfigure}		
	\caption{Convergence study of the error measured in the $L^2$, $H^1$ and $H^2$ norms in the non-matching case for nine patches example for different B-splines of degree $p=2,3$. Comparison of a classic penalty method, the scaled version with respect to the problem parameters proposed in~\citep{Herrema2019} (\textit{scaled}) and our projection approach (\textit{proj}).}
\label{fig:convergence_bilaplace_curved_9patches}
\end{figure}

In this example we consider the computational domain $\Omega = [0,3] \times [0,3]$ depicted in~\Cref{fig:setup_9patches}, divided into nine subdomains $\Omega^i$. 
Similarly to the previous example, all meshes are non-conforming at every coupling interface, where again an irrational factor of $\sqrt{2}/100$ has been used to shift the interface knots. 
The body source and boundary data are derived from the following analytical exact solution:
\begin{align*}
u^{\text{ex}} = \sin(\pi x) \cos(\pi x)  \, .
\end{align*}
Further, we set the Young's modulus to $E = 10^6 \, [Pa]$, the thickness of the plate to $t = 0.01 \, [m]$ and the Poisson's ratio to $\nu = 0 \, [-]$.
The convergence results of the error measured in the $L^2$, $H^1$ and $H^2$ are presented in~\Cref{fig:convergence_bilaplace_curved_9patches}, for splines of degree $p=2,3$. In this example we test the robustness of the method with respect to:
\begin{enumerate}
\item[$\bullet$] floating patches;
\item[$\bullet$] the presence of multiple cross-points where a constraint must be applied.
\end{enumerate} 
We observe again the expected asymptotic convergence rates of the error for all norms, where we remark that the method behaves optimally also for very coarse meshes, where locking phenomena are avoided.
Indeed, on one hand, we notice again that a classical ``vanilla'' choice of the penalty parameters yield a severe overconstraint of the solution space, resulting in a loss of accuracy of several order of magnitudes compared to the projection method. On the other hand, the scaling studied in~\citep{Herrema2019} leads to better results especially in the energy norm. However, for coarse meshes, we note that the method still suffers from locking, thus hindering the accuracy achievable by B-splines.

\subsection{A three patches example with a geometrically non-conforming interface}

\begin{figure}
	\centering
	\begin{subfigure}[t]{0.495\textwidth}
		\centering
	\includegraphics[width=0.75\textwidth]{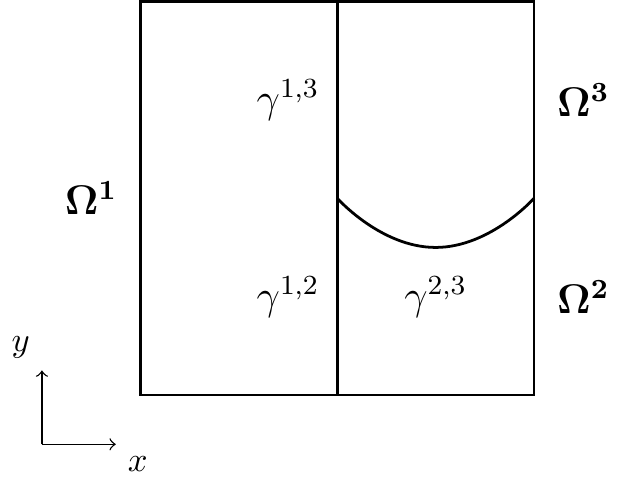}
		\caption{Initial subdomains and corresponding interfaces.}\label{fig:setup_3patches_initial}
	\end{subfigure}
	\hfill
	\begin{subfigure}[t]{0.495\textwidth}
		\centering
	\includegraphics[width=0.75\textwidth]{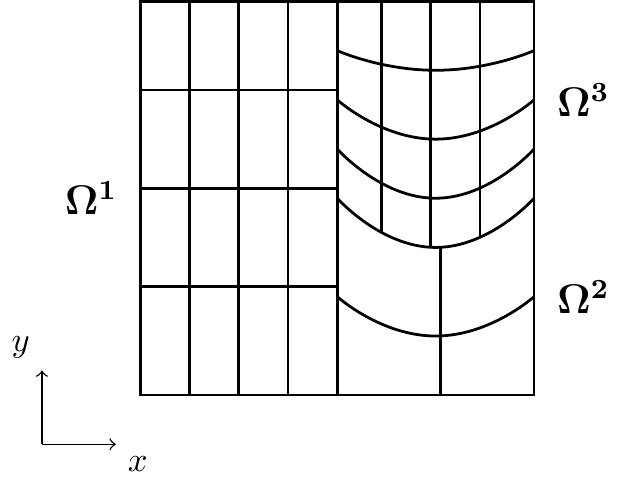}
		\caption{Non-conforming discretization.}\label{fig:setup_3patches_mesh}
	\end{subfigure}	
	\caption{Initial configuration and non-conforming discretization for the three patches example.}
\label{fig:setup_3patches}
\end{figure}

In this example we consider the computational domain $\Omega = [0,2] \times [0,2]$, split into three subdomains $\Omega^i$, see~\Cref{fig:setup_3patches_initial}. The initial non-conforming discretization used in the following is depicted in~\Cref{fig:setup_3patches_mesh},  where the interface knots are again shifted by a factor of $\sqrt{2}/100$ to induce the non-conformity.
The peculiarity of this example is the presence of a geometrically non-conforming interface between the patches, which is further used to assess the robustness of our method.
\begin{remark}
Similarly to~\citep{Brivadis2015}, we define an interface as geometrically conforming if the pull-back with respect to both slave and master domains is an entire edge of each parametric domain $\widehat{\Omega}^i$. 
\end{remark}
\noindent Similarly to the previous examples, the exact solution reads:
\begin{align*}
u^{\text{ex}} = \sin(\pi x) \cos(\pi x)  \, ,
\end{align*}
from which the applied body load and imposed boundary conditions are derived.
Regarding the problem parameters, we set the Young's modulus to $E = 10^6 \, [Pa]$, the thickness of the plate to $t = 0.01 \, [m]$ and the Poisson's ratio to $\nu = 0 \, [-]$.

\begin{figure}
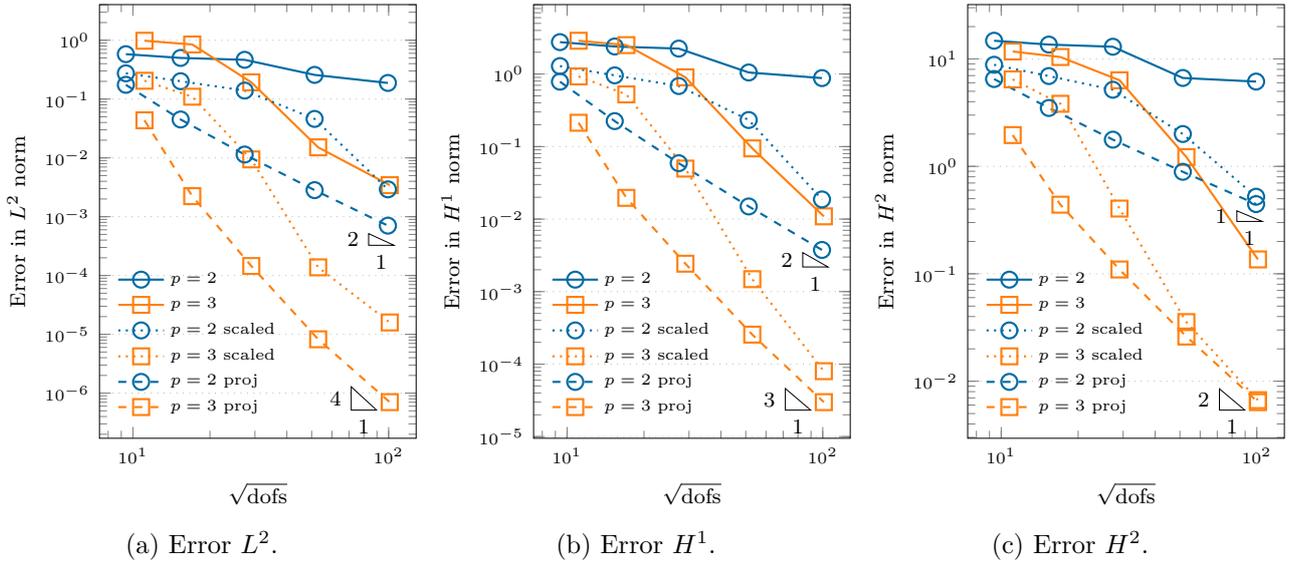

	\centering
	\begin{subfigure}[t]{0.32\textwidth}
		\centering
	\input{\graphDir/convergence_bilaplace_3patches_nonconforming_l2.tex}
		\caption{Error $L^2$.}
	\end{subfigure}
	\hfill
	\begin{subfigure}[t]{0.32\textwidth}
		\centering
	\input{\graphDir/convergence_bilaplace_3patches_nonconforming_h1.tex}
		\caption{Error $H^1$.}
	\end{subfigure}
	\hfill
	\begin{subfigure}[t]{0.32\textwidth}
		\centering
	\input{\graphDir/convergence_bilaplace_3patches_nonconforming_h2.tex}
		\caption{Error $H^2$.}
	\end{subfigure}		
	\caption{Convergence study of the error measured in the $L^2$, $H^1$ and $H^2$ norms in the non-matching case for the three patches example, B-splines of degree $p=2,3$. Comparison of a classic penalty method, the scaled version with respect to the problem parameters proposed in~\citep{Herrema2019} (\textit{scaled}) and our projection approach (\textit{proj}).}
\label{fig:convergence_bilaplace_3patches_nonmatching}
\end{figure}

The convergence results of the error measured in the $L^2$, $H^1$ and $H^2$ are presented in~\Cref{fig:convergence_bilaplace_3patches_nonmatching}, for splines of degree $p=2,3$. Analogously to our previous results, our method attains optimal rates of convergence, even in the presence of a geometrically non-conforming interface. Moreover, this numerical experiment confirms again that our method is insensitive to locking, starting from very coarse discretizations, where a substantial gain in accuracy per degree-of-freedom is observed.

\subsection{A flat L-bracket}

\begin{figure}
	\centering
	\begin{subfigure}[t]{0.495\textwidth}
	\centering
		\includegraphics[width=0.75\textwidth]{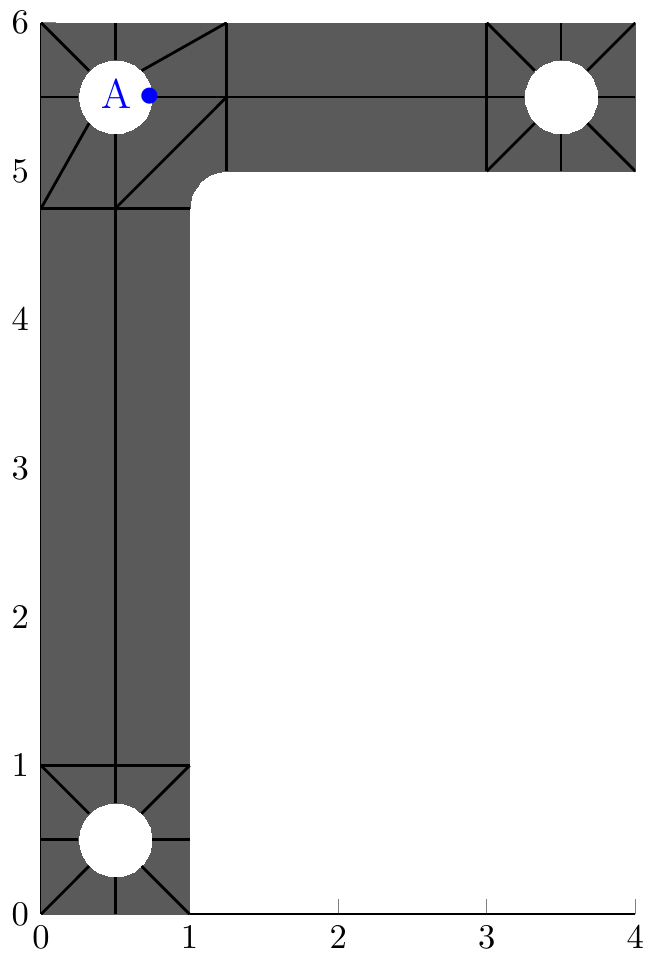}
		\caption{Multi-patch design.}\label{fig:setup_Lbracket_geo}
	\end{subfigure}
	\hfill
	\begin{subfigure}[t]{0.495\textwidth}
	\centering
		\includegraphics[width=0.75\textwidth]{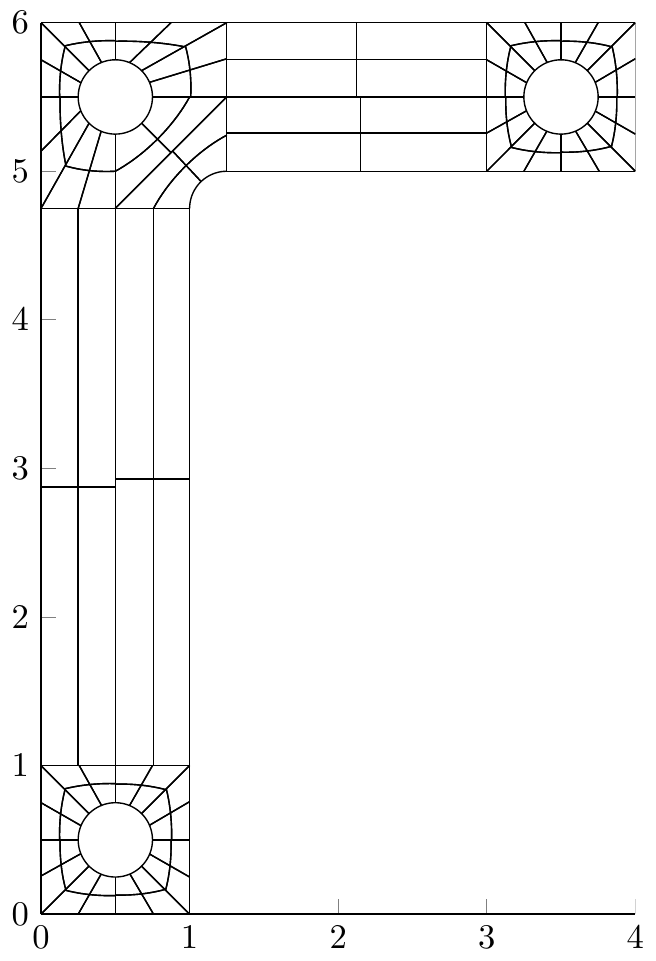}
		\caption{Non-conforming mesh.}
	\end{subfigure}
	\caption{Geometry setup and non-conforming discretization for the flat L-bracket example.}\label{fig:setup_Lbracket}
\end{figure}

\begin{figure}
	\centering
	\begin{subfigure}[t]{0.495\textwidth}
		\centering
	\includegraphics[width=0.75\textwidth]{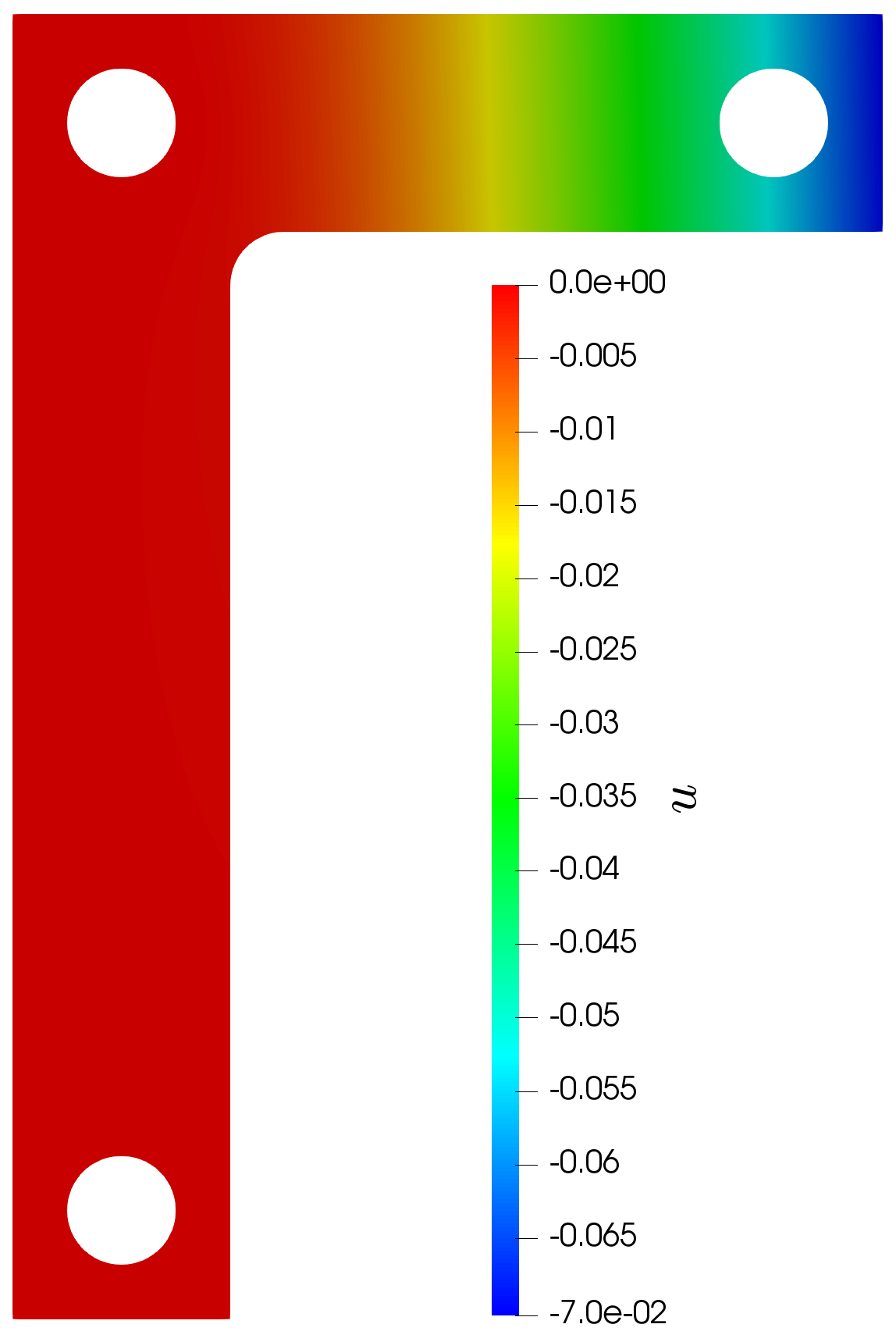}
		\caption{$p=2$.}
	\end{subfigure}
	\hfill
	\begin{subfigure}[t]{0.495\textwidth}
		\centering
	\includegraphics[width=0.75\textwidth]{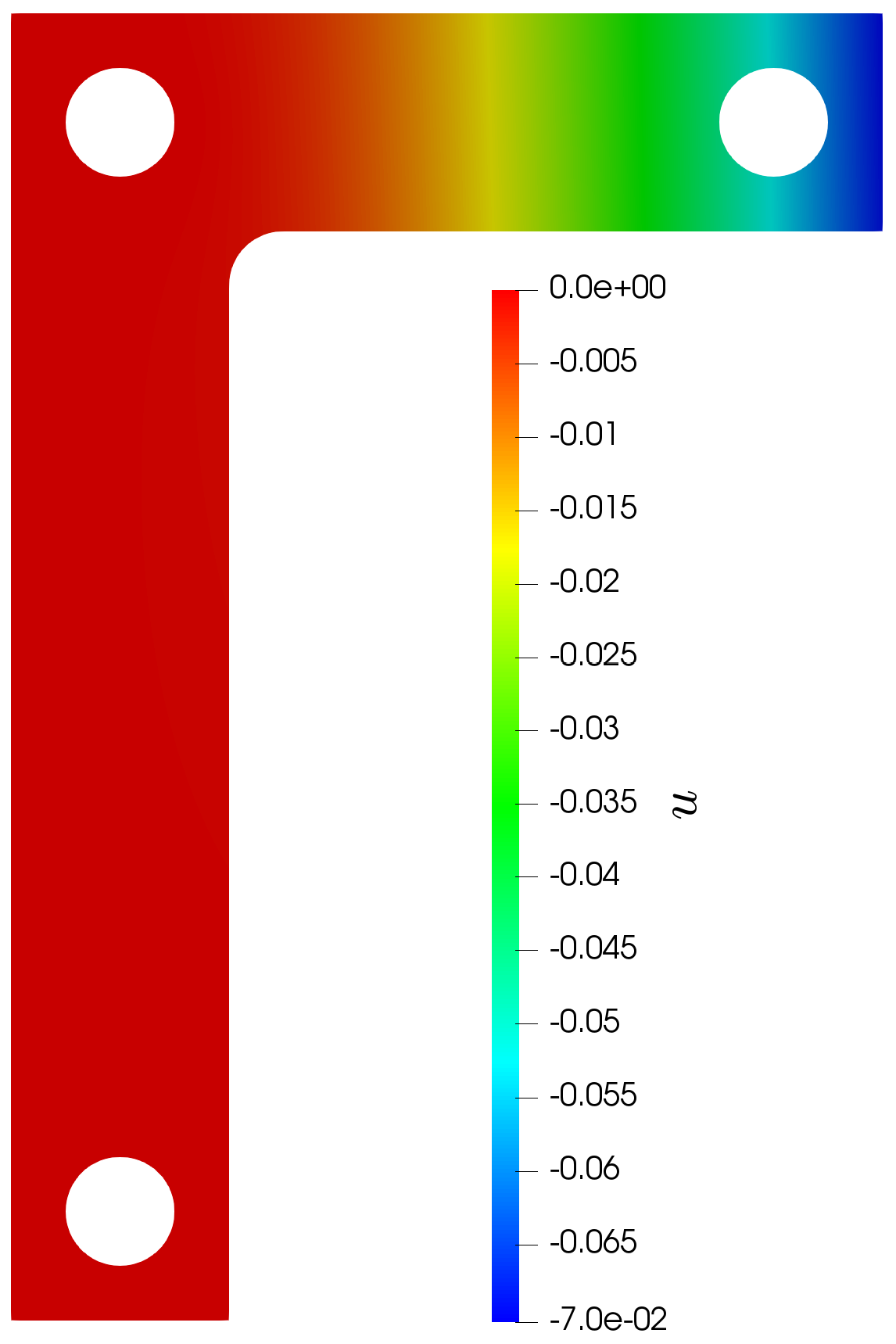}
		\caption{$p=3$.}
	\end{subfigure}	
	\caption{Solution contour for the flat L-bracket example, B-splines of degree $p=2,3$.}
\label{fig:sol_Lbracket}
\end{figure}

\begin{figure}
	\centering
	\begin{subfigure}[t]{0.325\textwidth}
		\centering
	\includegraphics[width=0.985\textwidth]{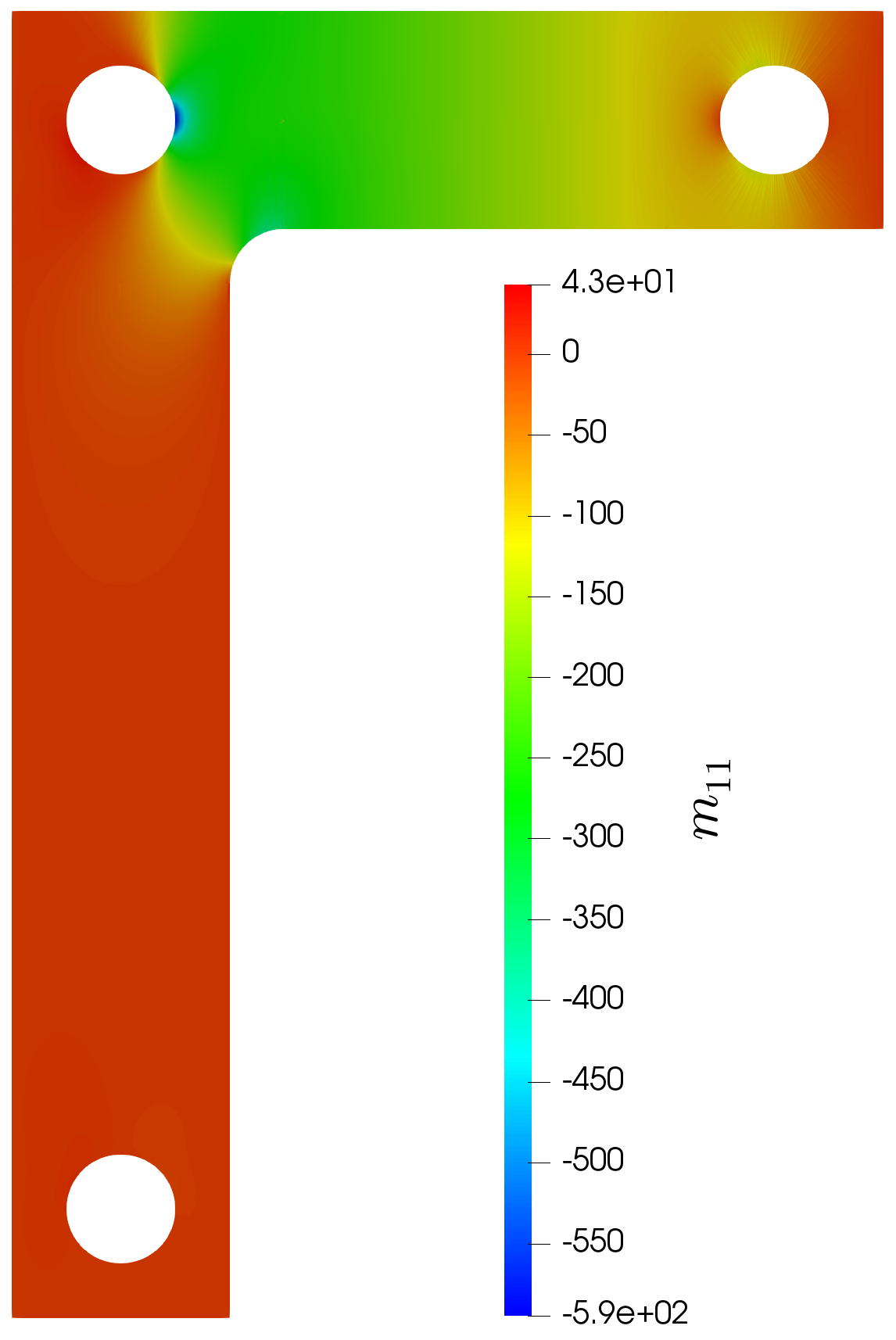}
		\caption{$m_{11}, p=2$.}
	\end{subfigure}
	\hfill
	\begin{subfigure}[t]{0.325\textwidth}
		\centering
	\includegraphics[width=0.985\textwidth]{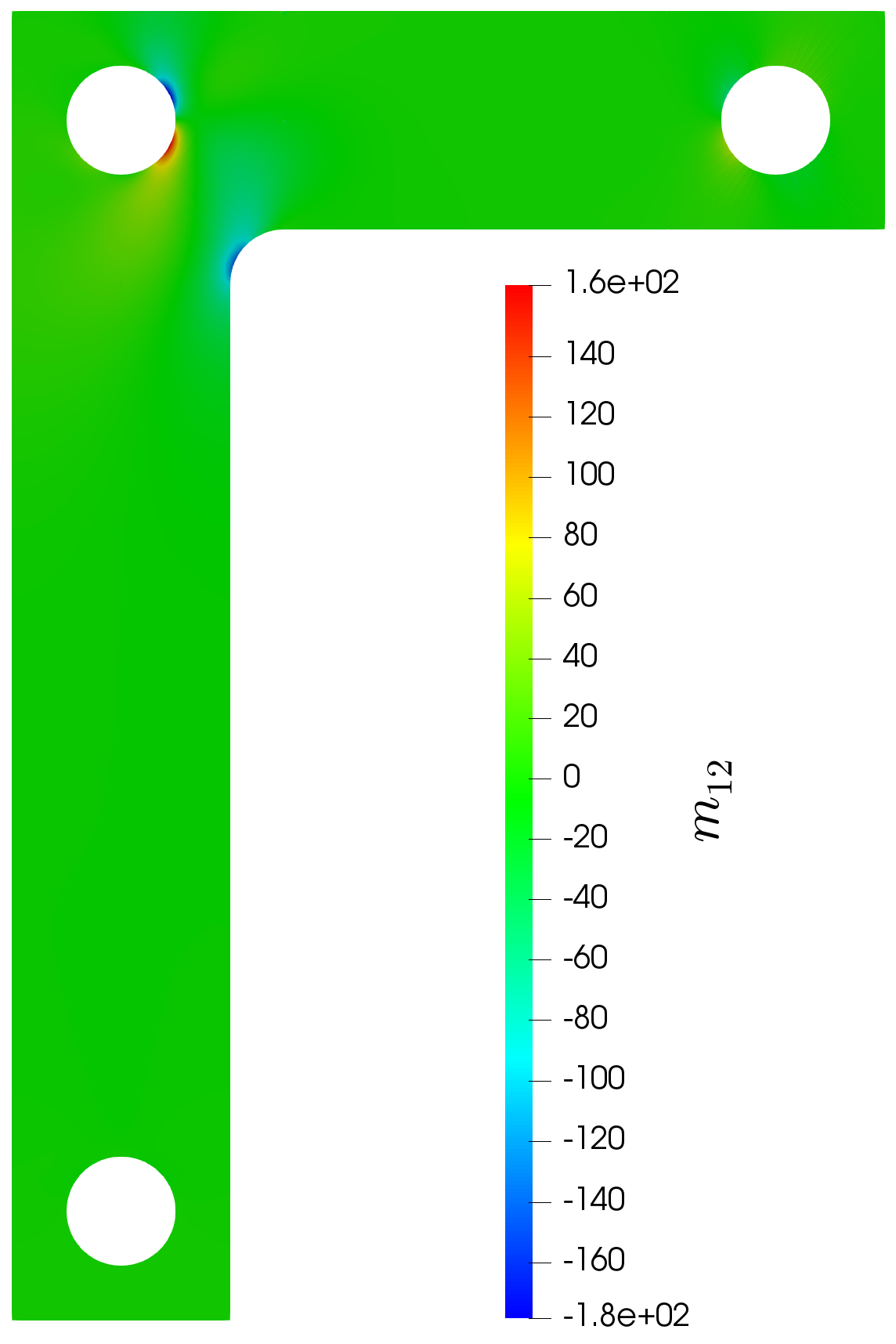}
		\caption{$m_{12}, p=2$.}
	\end{subfigure}	
	\begin{subfigure}[t]{0.325\textwidth}
		\centering
	\includegraphics[width=0.985\textwidth]{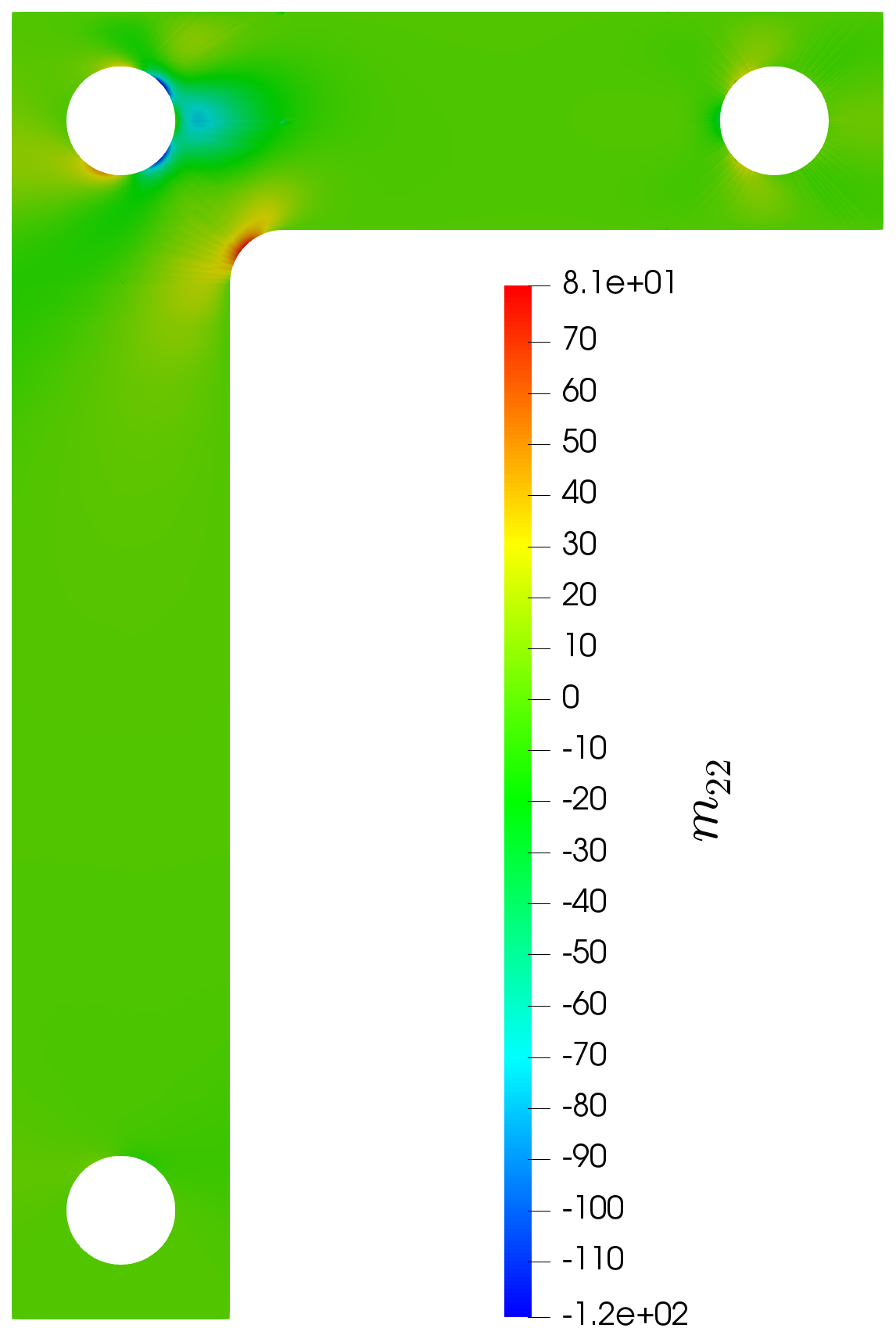}
		\caption{$m_{22}, p=2$.}
	\end{subfigure}
	\hfill
	\begin{subfigure}[t]{0.325\textwidth}
		\centering
	\includegraphics[width=0.985\textwidth]{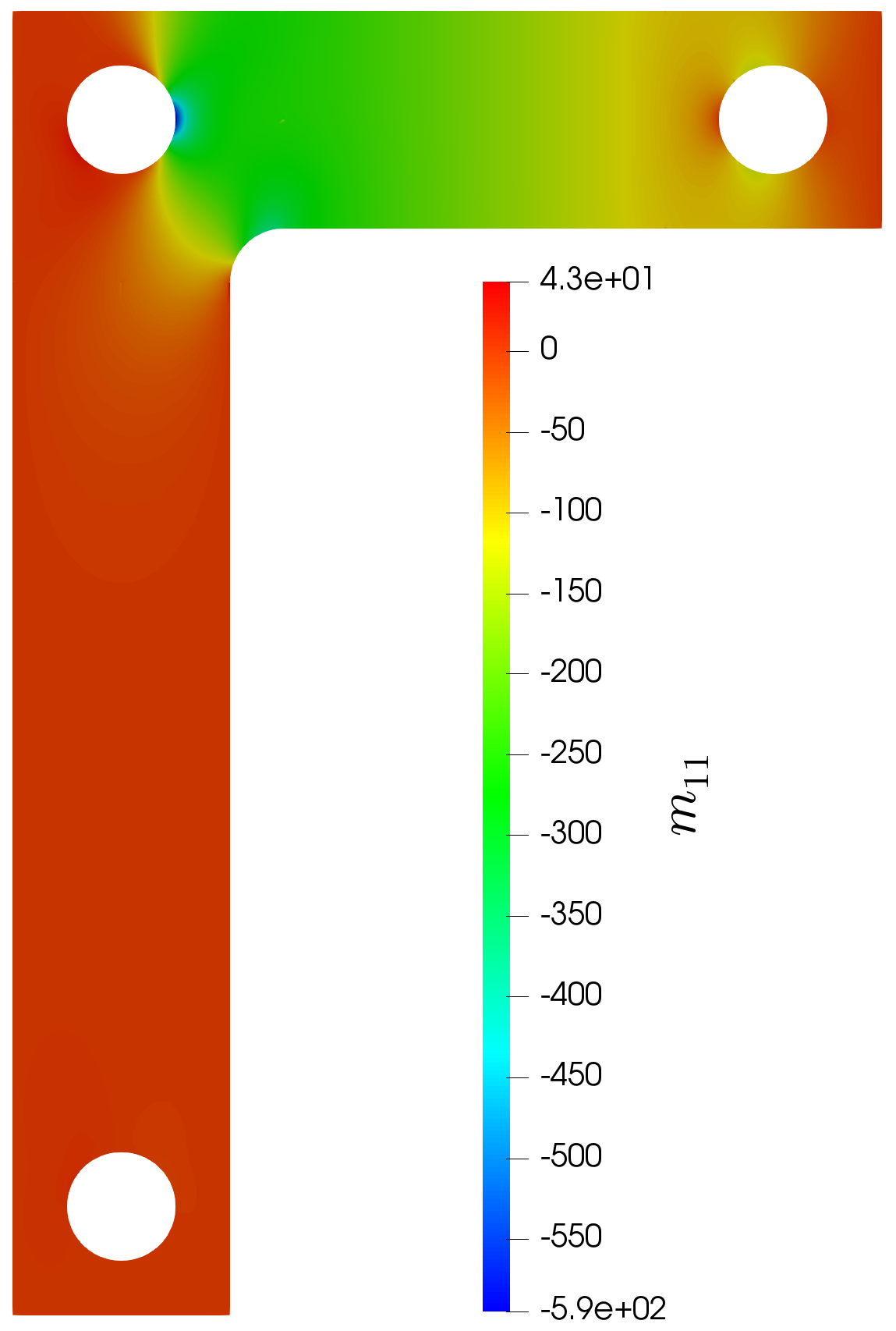}
		\caption{$m_{11}, p=3$.}
	\end{subfigure}	
		\begin{subfigure}[t]{0.325\textwidth}
		\centering
	\includegraphics[width=0.985\textwidth]{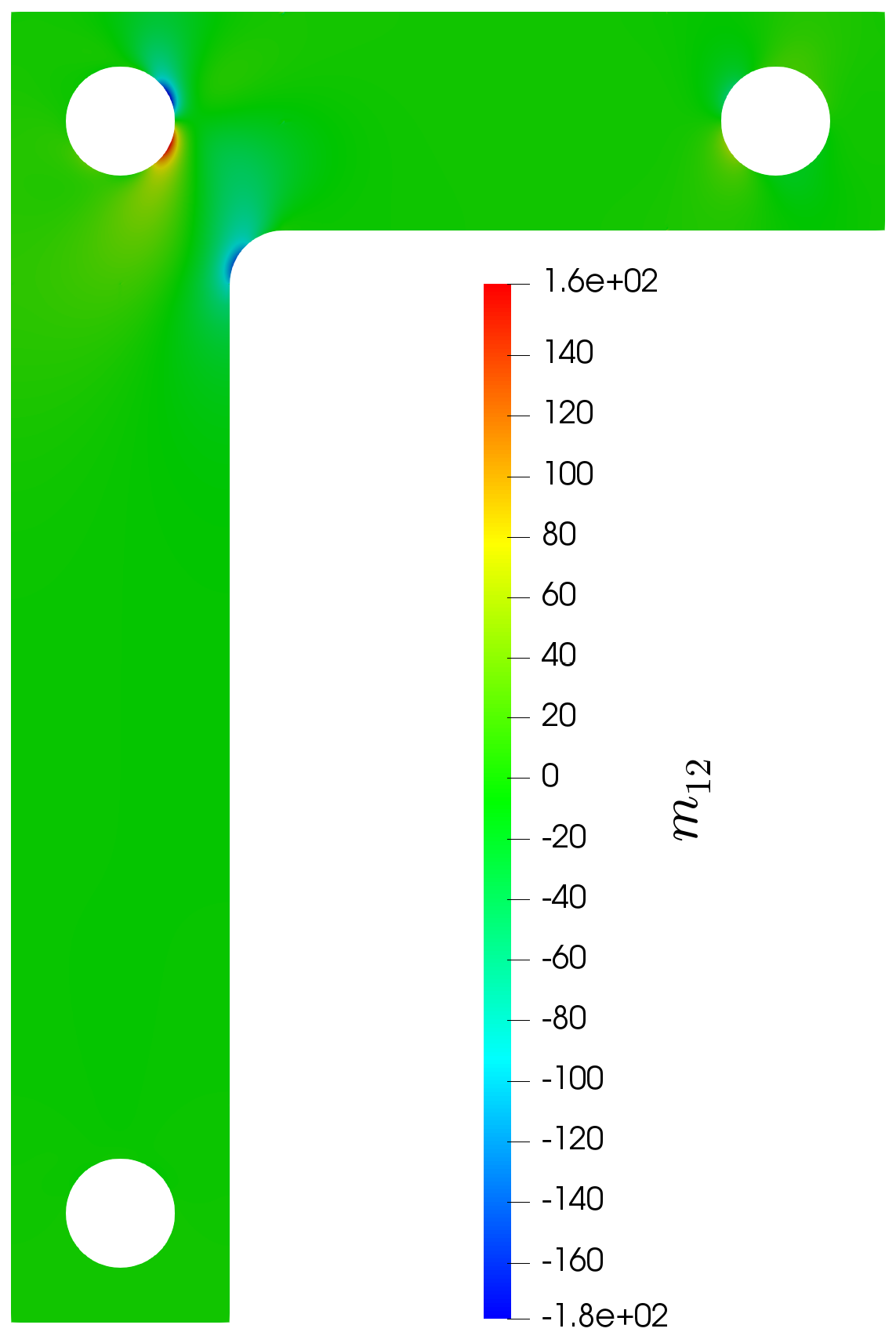}
		\caption{$m_{12}, p=3$.}
	\end{subfigure}
	\hfill
	\begin{subfigure}[t]{0.325\textwidth}
		\centering
	\includegraphics[width=0.985\textwidth]{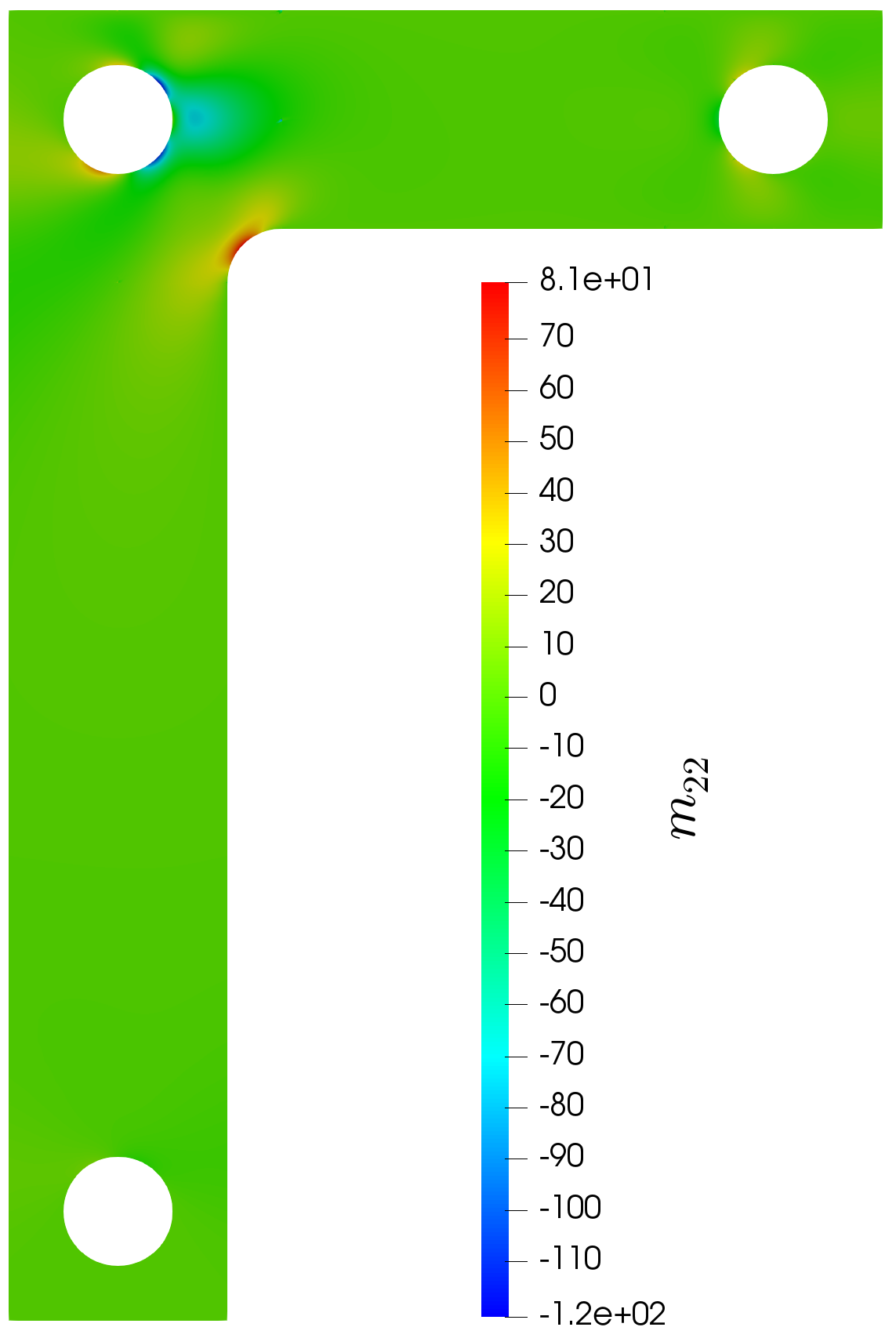}
		\caption{$m_{22}, p=3$.}
	\end{subfigure}		
	\caption{Components of the bending stress tensor $\mathbf{m}$ for the flat L-bracket example, B-splines of degree $p=2,3$.}
\label{fig:sol_stress_Lbracket}
\end{figure}

\begin{figure}
	\centering
	\input{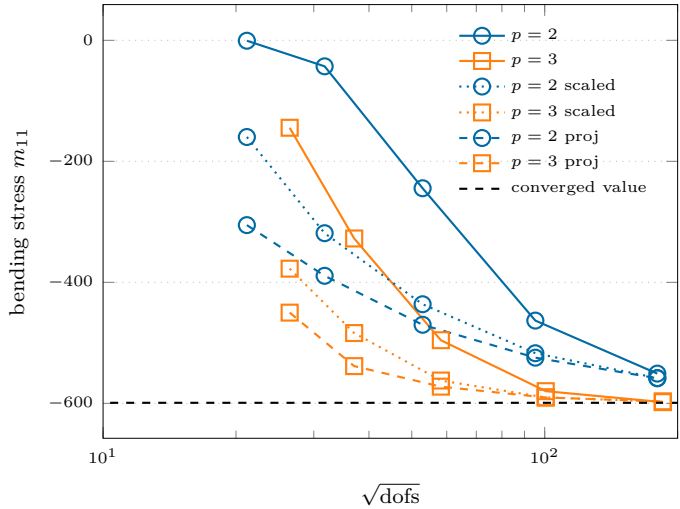}	
	\caption{Convergence study of the stress component $m_{11}$, evaluated at point A in~\Cref{fig:setup_3patches_initial}, for the flat L-bracket example for different B-splines of degree $p=2,3$. Comparison of a classic penalty method, the scaled version with respect to the problem parameters proposed in~\citep{Herrema2019} (\textit{scaled}) and our projection approach (\textit{proj}).}
\label{fig:convergence_m11}
\end{figure}

The last example we present is meant to show the applicability of the method to more complex multi-patch geometries. Analogously to the example studied in~\citep{Benzanken2017}, we modeled a flat L-bracket with 28 patches, coupled along 34 interfaces, as depicted in~\Cref{fig:setup_Lbracket}.
We applying a constant line load of $100 \, [N/m]$ in the negative $z$-direction on the upper right edge and we impose clamped boundary conditions on the entire boundary of the upper left and lower left holes, respectively. 
Further, we set the Young's modulus to $E = 200 \cdot 10^9 \, [Pa]$, the thickness of the plate to $t = 0.01 \, [m]$ and the Poisson's ratio to $\nu = 0 \, [-]$.
The solution field obtained with B-splines of degree $p=2,3$ is depicted in~\Cref{fig:sol_Lbracket}, where we remark the smoothness of the obtained solution, especially across the coupling interfaces.
In~\Cref{fig:sol_stress_Lbracket} we also plot the bending stress tensor $\mathbf{m}$, where its components are defined as:
\begin{align*}
m_{ij} = D \left( \nu \delta_{ij} u_{kk} + (1 - \nu) u_{ij} \right) \, ,
\end{align*}
and where $\delta_{ij}$ denotes the standard Kronecker delta. We obtain again a smooth stress field, where no visible spurious oscillations are introduced by the proposed coupling strategy.
Finally, in~\Cref{fig:convergence_m11}, we plot the convergence results of the stress component $m_{11}$, evaluated at point A marked in~\Cref{fig:setup_Lbracket_geo}, as a function of the number of dofs on a series of uniformly refined meshes. 
We note that for the classical penalty approach, and only for this example, we have tuned the penalty parameters to converge towards the reference value, where we have set $\alpha_{\text{defl}}^\ell = 10^4 E \, , \alpha_{\text{rot}}^\ell = E, \ell = 1,\ldots,L$.
This example highlights once again the gain in accuracy achieved on coarse meshes by the proposed method, also for point-wise quantities of interest.

\section{Conclusions} \label{sec:conclusions}

In this work we have introduced a simple methodology for the $C^1$-coupling of isogeometric patches based on the $L^2$-projection of suitable super-penalty terms in the context of Kirchhoff plates. The method does not suffer from locking phenomena, even in the case of severe non-matching discretization, where optimal rates of convergence of the error measured in the $L^2$, $H^1$ and $H^2$ norms have been attained also on very coarse meshes and a substantial gain in accuracy per degree-of-freedom has been observed compared to a classical penalty approach and to the scaled choice of parameters presented in~\citep{Herrema2019} in the scope of Kirchhoff-Love shells. The method turns out to be particularly effective for moderate spline degrees $p=2,3$. 
Our choice of parameters is completely determined by the problem definition and is based upon the underlying perturbed saddle point formulation associated to the plate, from which the two Lagrange multipliers are eliminated and the magnitude of the corresponding perturbations gives us insights on how to appropriately select the penalty factors. 
Then, to mitigate the detrimental effects of this choice on the condition number of the system matrix, we have combined the nested block preconditioner introduced in~\citep{Liu2019} with a preconditioner based on the Fast Diagonalization algorithm tailored for isogeometric Kirchhoff plates, inspired by the strategy in~\citep{Tani2016}. 

To conclude, we have demonstrated numerically the applicability and robustness of the proposed projected super-penalty approach for coupling Kirchhoff plates discretized by non-conforming isogeometric patches, where the method does not show any locking also on very coarse meshes.

\section*{Acknowledgements} 
The authors would like to thank Prof.~Alessandro Reali and Dr.~Rafael V\'{a}zquez for the fruitful discussions on the subject of this paper and Luca Pegolotti for his help with the implementational aspects of the preconditioner. The authors L.~Coradello and A.~Buffa gratefully acknowledge the support of the European Research Council, via the ERC AdG project CHANGE n.694515. The author A.~Buffa also gratefully acknowledges the support of the H2020-FetOpen-Ria project ADAM$^2$ n.862025.

\bibliographystyle{plainnatnourl}
\bibliography{library.bib}

\end{document}